\documentclass[a4paper,12pt]  {article}
\usepackage{amsmath,amsfonts}
\usepackage{amsthm}
\usepackage{amssymb}
\usepackage{graphicx}
\usepackage{verbatim}

\theoremstyle{plain}
\newtheorem{thm}{Theorem} 
\newtheorem*{thmA}{Theorem A} 
\newtheorem*{thmB}{Theorem B} 
\newtheorem*{thmC}{Theorem C} 
\newtheorem*{thmD}{Theorem D} 
\newtheorem*{thmE}{Theorem E} 

\newtheorem*{conj}{Conjecture}
\newtheorem*{cornn}{Corollary}
\newtheorem{lem}{Lemma}
\newtheorem{prop}{Proposition}

\theoremstyle{definition}
\newtheorem*{defn}{Definition}
\newtheorem*{defns}{Definitions}
\newtheorem*{rem}{Remark}
\newtheorem*{rems}{Remarks}

\newtheorem{exnum}{Example}

\begin{document}

\pagestyle{myheadings}
\markboth{Real algebraic curves and operator theory}
{Real algebraic curves and operator theory}

%
\title{Real Separated Algebraic Curves, \\
Quadrature Domains, Ahlfors Type Functions \\
and Operator Theory}
\author{Dmitry V. Yakubovich${}^1$ \\
Departamento de Matem\'aticas \\
Universidad Aut\'onoma de Madrid \\
Cantoblanco 28049, Madrid, Spain \\
dmitry.yakubovich@uam.es
}


\maketitle

\

\

\

\footnotetext[1] {This work was finished under the support of the
Ram\'on and Cajal Programme (2002) by the Ministry of the Science
and Technology of Spain.}


\newpage

\

\

\vskip4cm

\begin{abstract}
The aim of this paper is to inter-relate several
algebraic and analytic objects, such as
real-type algebraic curves, quadrature domains, functions
on them and rational matrix functions with special properties, and
some objects from Operator Theory, such as vector Toeplitz operators
and subnormal operators. Our tools come from operator theory, but
some of our results have purely algebraic formulation.
We make use of Xia's theory of subnormal operators and
of the previous results by the author in this direction.
We also correct (in Section 5) some inaccuracies in the works
of \cite{Ysubone}, \cite{Ysubtwo} by the author.
\end{abstract}

\


\

{Primary 30F50; Secondary 14H45, 14P05, 14M12, 46N99, 47B20, 65E05}

\newpage

%


\newcommand{\nc}{\newcommand}
\newcommand{\rnc}{\renewcommand}

\nc{\tht}{\theta} \nc{\clos}{\operatorname{clos}}
\nc{\wind}{\operatorname{wind}} \nc{\supp}{\operatorname{supp}}
\nc{\Ree}{\operatorname{Re}} \nc{\Res}{\operatorname{Res}}
\nc{\spaa}{{\operatornamewithlimits{span}}}
\nc{\spa}{\overline{\operatornamewithlimits{span}}}
\nc{\sm}{{\setminus}} \nc{\dist}{\operatorname{dist}}
\nc{\Ran}{\operatorname{Range}} \nc{\rank}{\operatorname{rank}}
\nc{\inter}{\operatorname{int}} \nc{\Range}{\operatorname{Range}}
\nc{\Pl}{\operatorname{Pol}} \nc{\Ker}{\operatorname{Ker}}
\nc{\Hol}{\operatorname{Hol}}
\nc{\Codim}{\operatorname{Codim}}
\nc{\const}{\operatorname{const}} \rnc{\phi}{\varphi}
\rnc{\epsilon}{\varepsilon} \nc{\eps}{\varepsilon}
\nc{\si}{\sigma} \nc{\Pzt}{P_{z,\,t}}

\nc{\defin} {{\,\overset {\text {\rm def} }{=}\,}}

\nc{\smma}{\smallmatrix} \nc{\esmma}{\endsmallmatrix}
\nc{\gews}{\geq} \nc{\lews}{\leq} \nc{\al}{\alpha} \nc{\be}{\beta}
\nc{\ga}{\gamma}
\nc{\Ga}{\Gamma} \nc{\de}{\delta} \nc{\ze}{\zeta}
\nc{\pl}{{\partial}} \nc{\vardelta}{{\bold\delta}}
\nc{\Del}{\Delta} \nc{\De}{\Delta} \nc{\wDe}{\widehat\Delta}
\nc{\wDeR}{\wDe_{\mathbb R}} \nc{\wDedeg}{\wDe_{\rm deg}}
\nc{\wDendeg}{\wDe_{\rm ndeg}} \nc{\wDendegp}{\wDe_{\rm ndeg}^+}
\nc{\Dendeg}{\De_{\rm ndeg}} \nc{\BP}{{\mathbb P}}
\nc{\BR}{{\mathbb R}} \nc{\vph}{\vphantom{\hat{\hat L}}}
\nc{\vvph}{\vphantom{\hat L}}

\nc{\PHd}{P_{H'}}

\nc{\gash}{\partial\wDe_+} \nc{\Bsh}{B^\#} \nc{\CA}{\mathcal A}
\nc{\cD}{\mathcal D} \nc{\CU}{\mathcal U} \nc{\CV}{\mathcal V}
\nc{\CL}{\mathcal L} \nc{\CK}{\mathcal K} \nc{\CN}{\mathcal N}
\nc{\CP}{\mathcal P} \nc{\CM}{\mathcal M} \nc{\CW}{\mathcal W}
\nc{\E}{\mathcal E} \nc{\F}{\mathcal F} \nc{\K}{\mathcal K}
\nc{\bK}{\overline{\mathcal K}} \nc{\BC}{\mathbb C}
\nc{\BD}{\mathbb D} \nc{\BN}{\mathbb N} \nc{\BQ}{\mathbb Q}
\nc{\BT}{\mathbb T}
\nc{\BZ}{\mathbb Z} \nc{\gA}{\gothic A}
\nc{\wt}{\widetilde } \nc{\wtH}{\widetilde H} \nc{\wtN}{\widetilde
N} \nc{\wtR}{\widetilde R} \nc{\wtS}{\widetilde S}
\nc{\wtU}{\widetilde U} \nc{\wtV}{\widetilde V} \nc{\V}{\mathcal
O} \nc{\wtW}{\widetilde W} \nc{\wtX}{\widetilde X}
\nc{\wtG}{\widetilde G} \nc{\wtBD}{\widetilde\BD}
\nc{\ovlR}{\overline{R}} \nc{\whR}{\widehat{R}} \nc{\CLM}{\mathcal
L(M)} \nc{\gasym}{\gamma_{\text{symm}}} \nc{\Om}{\Omega}
\nc{\wtOm}{\widetilde \Om} \nc{\om}{\omega}

\nc{\La}{\Lambda} \nc{\la}{\lambda} \nc{\sCLa}{\sigma_C(\La)}
\nc{\sLa}{\sigma(\La)} \nc{\sLas}{\sigma(\La^*)}
\nc{\Htwos}{{H^2_s}} \nc{\Htwok}{{H^2_k}} \nc{\Htwom}{{H^2_m}}
\nc{\Htwomin}{{H^2_{m,-}}} \nc{\Qfl}{Q^\flat} \nc{\Le}{L^2(e)}
\nc{\CLAS}{C(z-\La)^{-1}+\La^*}
\nc{\lw}{\leftrightarrow} \nc{\NDARN}{{\operatorname{NDARN}}}
\nc{\ind}{{\operatorname{ind}}} \nc{\NDRN}{{\operatorname{NDRN}}}
\nc{\NDARNmm}{{\operatorname{NDARN}_m}}
\nc{\NDRNmm}{{\operatorname{NDRN}_m}}
\nc{\NDARNm}{{\operatorname{NDARN}_m}}
\nc{\NDRNm}{{\operatorname{NDRN}_m}}
\nc{\codim}{{\operatorname{codim}}}
\nc{\Area}{{\operatorname{Area}}}
\nc{\trace}{{\operatorname{trace}}} \nc{\Hinf}{{H^\infty}}
\nc{\Hinfmm}{{H^\infty_{m\times m}}}
\nc{\Hinfms}{{H^\infty_{m\times s}}}
\nc{\Hinfss}{{H^\infty_{s\times s}}}
\nc{\Hinfmk}{{H^\infty_{m\times k}}}
\nc{\Hinfkm}{{H^\infty_{k\times m}}}
\nc{\Linfmm}{{L^\infty_{m\times m}(\BT)}}

\nc{\partR}{{\partial R}}

\nc{\partBD}{\BT}

\nc{\partwDep}{\partial\wDe_+} \nc{\alj}{{\al_j}} \nc{\kj}{{k_j}}


\nc\wh {\widehat} \nc\wH {\wh H} \nc{\smc}{\textsc}

\nc\beqn{\begin{equation}} \nc\neqn{\end{equation}}
\nc{\beqnay}{\begin{eqnarray}}   \nc{\neqnay}{\end{eqnarray}}
\nc{\beqnays}{\begin{eqnarray*}} \nc{\neqnays}{\end{eqnarray*}}
\nc{\barr}{\begin{array}}        \nc{\narr}{\end{array}}
\nc\nn{\nonumber}


\setcounter{section}{-1}

\section{Introduction}

In this paper, we discuss the inter-relations
between the following objects:

(1) Separated real algebraic curves in $\BC^2$;

(2) Algebraic curves in $\BC^3$ of special type, which we will
call Ahlfors type curves;

(3) Quadrature domains;

(4) Rational matrix functions of a certain class;

(5) The corresponding (analytic) Toeplitz operators on vector
Hardy spaces $H^2_m$, $m\in\BN$;

(6) Subnormal operators with finite rank self-commutators and
isometries that commute with them.

Each of the objects we discuss depends on a finite number of real
parameters, and some of the
connections we speak about are given by explicit
formulas (see, for instance, \S7).

Here, in the Introduction, we
will define these objects and formulate some of our main
results. More background and more explanations will be given in
next sections.

A polynomial $Q(z,w)$ will be called {\it of real type} if
it has a form
\beqn
Q(z,w)=\sum_{j=0}^n\sum_{k=0}^n a_{jk}z^j w^k,
\label{AEE}
\neqn
where $n\ge1$ is an integer, $a_{jk}\in\BC$,  and     
\beqn
a_{kj}=\bar a_{jk}, \qquad 0\lews j,k\lews n.
\label{ajk}
\neqn

The linear invertible substitution $z=x+iy$, $w=x-iy$, $x,y\in\BC$
transforms each real-type polynomial into a real polynomial in variables
$x$, $y$ (the converse also is true).

Consider the algebraic curve
$$
\De=\bigl\{(z,w)\in\BC^2:Q(z,w)=0\bigr\}
\label{AF}
$$
that corresponds to $Q$.
Since $Q(\bar w,\bar z)=\overline{Q(z,w)}$, the anti-analytic
involution
\beqn
\de=(z,w)\mapsto\de^*\defin(\bar w,\bar z)
\label{invo}
\neqn
maps $\De$ onto
itself. We will call $\De$ {\it a real-type algebraic curve} if it is
given by an equation in the form (\ref{AEE}), subject to
(\ref{ajk}). Equivalently, an algebraic curve in $\BC^2$ is of
real type if
its equation is invariant under the involution (\ref{invo}).

The polynomial $Q$ admits a unique decomposition
$Q=\prod_{j=1}^T Q_j^{k_j}$ into a product of
irreducible polynomials \cite{FUL}; we will
call algebraic curves
$\De_j=\{Q_j=0\}$
in $\BC^2$
\textit{the irreducible pieces} of the curve $Q$
and write symbolically $\De=\De_1^{k_1}\cup\dots\cup\De_T^{k_T}$.
Each of $\De_j$ essentially is a compact (unbordered) Riemann surface
of finite genus. That is, there is only a finite number of singular points
on $\De_j$, where $\frac{\pl  Q_j}{\pl  z}=\frac{\pl  Q_j}{\pl  w}=0$,
and by picking out these points from $\De_j$ and adding a finite number
of new (``ideal'') points one gets a compact (unbordered)
abstract Riemann surface $\wDe_j$.
This procedure is unique, and this
surface is called a
\textit{desingularization} of $\De_j$ (see \cite{FUL}).
Coordinate functions $z$, $w$ are globally meromorphic functions on $\wDe_j$.

We assume the desingularizations $\wDe_1\dots\wDe_k$ to be disjoint.
The formula $\de\mapsto(z(\de),w(\de))$
defines a continuous ``projection'' of
$\wDe$ onto $\De$. The involution (\ref{invo}) is also defined on  $\wDe$.

We put $\wDeR=\big\{\de\in\wDe: \de=\de^*\big\}$ to be the real part of
$\wDe$. (In general, a point of a compact Riemann surface with anti-analytic involution
is called
\textit{real} if it is invariant under the involution.)

\begin{defns}
1) A component $\De_j$ of $\De$ is called \textit{degenerate}
if it has either the form
$(z-a)^{k_j}=0$ or  $(w-\bar a)^{k_j}=0$, and
\textit{non-degenerate} in all other cases.

We put $\wDe_{\text{deg}}$ to be the union of
all degenerate components of $\wDe$, and
$\wDe_{\text{ndeg}}$ to be the union of
all its non-degenerate components (with their multiplicities).

2) An irreducible
real-type algebraic curve $\wDe$ will be
called {\it separated}
if the real part $\wDeR$
separates $\wDe$ in the topological case.
We will say that a real-type algebraic curve
is separated if all its non-degenerate pieces
are separated.

3) Let  $\wDe$ be a separated real-type algebraic curve.
We call $\wDe$
\textit{pole definite} if no pole of
$z(\cdot)$ lies on the real part of $\wDe$ and for each
non-degenerate irreducible piece $\wDe_j$ of $\wDe$,
all poles of $z(\cdot)$ belong to
the same connected component of
$\wDe_j\sm\wDeR$.
\end{defns}

If $\wDe$ is separated, then the involution (\ref{invo}) maps
each of
its non-degenerate irreducible pieces onto itself.
The general theory \cite{NAT} of Klein surfaces
(Riemann surfaces with anti-analytic involution)
implies
that in this case, for any
piece $\wDe_j$ of $\wDendeg$,
the complement $\wDe_j\sm\wDeR$ has exactly two connected components.
In the case when $\wDe$ is pole definite, we will call these
connected components $\wDe_j^+$ and $\wDe_j^-$, assuming that
$z(\wDe_j^+)$ is bounded and $z(\wDe_j^-)$ is not.
The involution
interchanges $\wDe_j^+$ with $\wDe_j^-$, so that
they can be called ``halves'' of the piece
$\wDe_j$.

The coordinate function $\de\mapsto z(\de)$
is  bounded and analytic on $\wDe_j^+$.

We put
$$
\wDe_\pm=\bigcup \wDe_j^\pm,
$$
so that we have a disjoint union
$$
\wDendeg=\wDe_+\cup\wDe_-\cup\wDeR.
$$

\

{\bf Quadrature domains. } Suppose $\Om$ is a bounded domain in
$\BC$ and there are points $z_k$ in $\Om$ and complex constants
$c_{jk}$ such that an identity
\beqn
\iint_\Om
f\,dx\,dy=\sum_{k=1}^s\sum_{j=0}^{r_k-1}c_{jk}f^{(j)}(z_k)
\label{tAA}
\neqn
holds for all analytic functions $f$ in $L^1(\Om)$. Then $\Om$ is called
a quadrature domain. We will call points
$z_k$ \textit{the nodes}  and
the coefficients $c_{jk}$
\textit{the weights} of our domain $\Om$.

Quadrature domains possess many interesting and
intriguing properties. After the
pioneering work \cite{AhSh} by Aharonov and Shapiro,
in the last 20 - 30 years,
quadrature
domains have been related with such diverse fields as algebraic geometry
\cite{AhSh}, \cite{Shdou}--\cite{GP},
potential theory and different problems in fluid dynamics
\cite{Sak}, \cite{Cro3}--\cite{Cro}, \cite{Richsn},
moment problems
\cite{PUT}, \cite{PUTJFA}, \cite{GustHeMiPutinar},
extremal problems for univalent functions (studied by
Aharonov, Shapiro and Solynin, see \cite{AhShSoly}
and references therein). They also have close relations with
subnormal and hyponormal operators
``of finite type''
\cite{Xone}--\cite{Xrev}, \cite{GP},
\cite{PUT}, \cite{Ysubone}--\cite{Yqua}.  We refer to the
recent book \cite{EbenGuKhPutinar} for more information.

A function $w(z)$ on $\clos \Om$ is called
\textit{a Schwartz function} of $\Om$ if $w$ is holomorphic on
$\Om$, except for finitely many poles,
continuous on the boundary of $\Om$ and satisfies $w(z)=\bar z$,
$z\in\partial\Om$. It has been known since the work
 \cite{AhSh} that $\Om$ is a
quadrature domain if and only if it possesses a
a Schwartz function. Moreover, in this case
the nodes of $\Om$ coincide with the poles of
the function $w$. If the poles are simple,
only the weights $c_{0k}$ are
present in (\ref{tAA}), and they are
proportional to the residues of
$w(z)$ at points $z=z_k$.

The following result relates
algebraic curves of the above form with quadrature domains.

\begin{thmA}[{\rm Aharonov and Shapiro \cite{AhSh}, Gustafsson \cite{Shdou}}]
If $\De$ is an irreducible separated real-type algebraic curve
and the coordinate function $z$  is injective on
$\wDe_+$, then the image $z(\wDe_+)$ is a quadrature domain.
Each quadrature domain is formed in this way.
\end{thmA}

In the situation of this theorem,
the Schwartz function $w(z)$ on $\Om$ coincides with the coordinate $w$
on our curve. More precisely,
the Schwartz function is $w\big((z|\wDe_+)^{-1}\big)$.
The Riemann surface $\wDe$ can be constructed
from the domain $\Om$ as
follows. Take one more copy $\wt\Om$ of $\Om$ and endow it with the conformal
structure, provided by the function $\bar z$.
Then $\wDe$ is isomorphic as a compact Riemann surface to
the so-called \textit{Schottky double} of $\Om$, which is obtained
by welding $\Om$ and $\wt\Om$ together along $\pl\Om$.
We refer to \cite{Shdou} for more details.
The irreducibility of the curve $\wDe$ follows from the fact that the
Schottky double is always a connected topological space.

\

{\bf Rational matrix functions. }
Suppose $F$ is a continuous $m\times m$ matrix function on the unit circle $\BT$.
Then eigenvalues of $F(t)$ depend continuously on
the point $t\in\BT$, that is, there  are
continuous functions $\ze_1(\theta),\dots,\ze_m(\theta)$,
$\theta\in[0,2\pi]$ such that for each $\theta$,
$F(e^{i\theta})$ has eigenvalues $\ze_1(\theta),\dots,\ze_m(\theta)$,
counted with
their multiplicities. In general,
$\ze_j(0)\ne\ze_j(2\pi)$. For a point
$z_0\in\BC$
such that $\det(F(t)-z_0I)$ does not vanish for $t\in\BT$,
we define \textit{the winding number} of the matrix function
$F$ around $z_0$ as the sum of the increments of the
argument of $\ze_j(\cdot)-z_0$:
\beqn
\label{wind}
\wind_F(z_0)
=\frac1{2\pi}\sum_{j=1}^m
\De_{[0,2\pi]}\arg\big(\ze_j(\cdot)-z_0\big).
\neqn
The number $\wind_F(z_0)$  equals to the winding number of the scalar function
$\theta\mapsto\det\big(F(e^{i\theta})-z_0 I \big)$ around
the origin. Hence it does not depend on the choice of
continuous branches of eigenvalues $\ze_1,\dots,\ze_m$.

We need some special classes of symbols.

\begin{defns}
Let $F$ be a square matrix function
on the (open) unit disc $\BD$. Function $F$ will be called

- analytic, if $F$ is bounded and analytic on $\BD$;

- normal, if the matrix $F(t)$ is
defined and is normal for almost every $t\in\BT$;

- non-degenerate, if for any constant $c$ in $\BC$, the determinant
$$
\det(F(t)-cI)
$$
is not identically zero on $\BD$.

We denote by
$\NDRNmm$ the class of all non-degenerate rational normal matrix functions
of size $m\times m$, and by $\NDARNmm$ the class of all non-degenerate
\textit{ analytic } rational normal matrix functions
of size $m\times m$ (so that $\NDARNmm$ is a subclass of $\NDRNmm$).
\end{defns}

If $F$ is rational, then it is non-degenerate iff for any $c$,
the above determinant is not identically zero on $\BC$.
A scalar rational function on $\BD$ without poles on
$\clos \BD$ is non-degenerate iff it is not constant.

We say that a domain $\Om$ in $\BC$ is \textit{$p$-connected}
(or \textit{has connectivity} $p$) if
the homology group ${\cal H}^1(\Om,\mathbb Z)$ is isomorphic to
$\mathbb Z^p$. A simply connected domain has
connectivity $0$ and a domain with one hole has connectivity $1$.
Sometimes the use of this term is different.

Aharonov and Shapiro also proved in \cite{AhSh} the following result.
\begin{thmB}[Aharonov and Shapiro]
A simply connected domain $\Om$ is a quadrature domain if and only if there is a
(scalar) rational function $g$, which is analytic on the closed
unit disc, univalent in the open unit disc $\BD$ and satisfies
$g(\BD)=\Om$.
\end{thmB}

In the situation of this theorem, equations
$z=g(t)$, $w=\overline{g\vvph(\bar t^{-1})}$, $t\in\BD$ define implicitly
the Schwartz function on $\Om$. If $\{t_k\}$ are the poles
of $g$ on the Riemann sphere, then the nodes of $\Om$ are
exactly the points $g(\bar t_k^{-1})$.
Here we denote by $z$ both the meromorphic
function $z(\de)$ on $\wDe$ and the independent variable of the
$z$-plane; it should not confuse the reader.

Our first result extends Theorem B to
the multiply connected case.

\begin{thm}
A bounded domain $\Om$ in $\BC$ is a quadrature domain
if and only if there is a natural number $m$ and a function
$F\in\NDARN_m$ with continuous branches
$\ze_1(\theta),\dots,\ze_m(\theta)$ of eigenvalues of
$F(e^{i\theta})$, $\theta\in[0,2\pi]$ such that

{\bf (1)} $\pl\Om=\bigcup_j\ze_j([0,2\pi])$;

{\bf (2)} $\wind_F(z)=1$ for $z\in\Om$;

{\bf (3)}  $\wind_F(z)=0$ for $z\in\BC\sm\clos\Om$.

\noindent If $\Om$ has connectivity $p$, then one can find
a function $F$ with these properties in $\NDARN_{p+1}$.
\end{thm}

If $\Om$ and $F$ are related as in this theorem, we will say that
\textit{the matrix function $F$ generates the domain $\Om$}.

In fact, in \S1  we will associate with any function $F\in\NDARN_m$
an algebraic curve $\De_{(2)}(F)$ in $\BC^2$
and an algebraic curve $\De_{(3)}(F)$ in $\BC^3$.
Theorem~3 and its Corollary assert
that an algebraic curve $\De$ in $\BC^2$ is
admissible, pole definite and separated iff
$\De=\De_{(2)}(F)$ for some
$F\in\NDARN_m$ (for a certain $m$).
Theorem 1 can be considered as a particular case
of this result. For the reader's convenience, we will prove
Theorem 1 independently of Theorem 3.

Now let us pass to the operator theory objects we will need.

\

{\bf Analytic vector Toeplitz operators.}
Let $m\gews 1$ be an integer. The \textit{vector Hardy space} is
$$
H^2_m\defin
\bigg\{\,
f(t)=\sum_{n\ge0}a_nt^n:\qquad a_n\in\BC^m,\;
\|f\|^2\defin\sum_{n\ge0}\|a_n\|^2<\infty\,
\bigg\};
$$
it is a Hilbert space of $\BC^m$-valued analytic functions
in the unit disc $\BD$. A function $f(t)$ in $H^2_m$ can also
be considered as a $\BC^m$-valued  $L^2$ function on the unit circle;
its values on the circle are radial limits of its values on $\BD$ a.e.
In this interpretation, the space $H^2_m$ becomes a closed subspace
of $L^2(\BC^m)$,
and in the above formula
$a_n$ are the Fourier coefficients of $f$.

 For $m=1$, $H^2_m$ is a classical
scalar-valued Hardy space $H^2$; in general, $H^2_m=\bigoplus_{j=1}^m H^2$.
We refer to \cite{DUR} for basic properties of Hardy spaces $H^p$.

The class $\Hinf$ consists of all bounded analytic functions in $\BD$; it is
equipped with the supremum norm. We denote by $\Hinfmk$ the class of
$m\times k$ matrix functions on the unit disc, whose entries are in $\Hinf$.
These matrix functions, certainly, also have boundary limit values a.e. on $\BT$.

Let $F\in\Hinfmm$. The analytic vector Toeplitz operator on the
vector Hardy space $H^2_m$
with the symbol $F$ is in fact a multiplication operator, which
acts by the formula
$$
T_Fg(t)=F(t)g(t), \qquad g\in H^2_m
$$
(the general definition of a vector Toeplitz operator will be given in \S3).

\

{\bf Subnormal operators. }
Let $H$ be a Hilbert space.
Throughout the article, we will deal only with separable complex Hilbert
spaces and bounded linear operators.
We denote by $\CL(H_1, H_2)$ the set of linear operators
acting from $H_1$ to $H_2$ and write $\CL(H)$ instead of
$\CL(H, H)$.

\begin{defns}
A linear operator $S$ acting on a
Hilbert space $H$ is called \textit{subnormal}
if there exist a larger Hilbert space $K$, $K\supset H$ and a normal operator
$N$ in $\CL(K)$ such that $NH\subset H$ and $S=N|H$.
In this case, we call $N$
\textit{a normal extension} of $S$.
We will say that $S$ \textit{has no point masses} if
it has  a normal extension $N$ that has no non-zero eigenvectors.
We call $S$ \textit{pure} if it has no nonzero invariant subspace, on
which it is normal.

We will say that a subnormal operator $S$
\textit{ is of finite type}
if it is pure and its self-commutator
$[S^*,S]\defin S^*S-SS^*$ has finite rank.
\end{defns}

Subnormal operators have been much investigated; we refer to
the book \cite{Conw} for a background.


In a general setting,
a kind of the spectral theory of subnormal operators
was developed by Xia in \cite{Xone}--\cite{Xadd}. In
author's previous work \cite{Ysubone}, \cite{Ysubtwo}, an
alternative exposition of Xia's theory
for subnormal operators of finite type was given.  A strong
two-sided relationship between
subnormal operators of finite type and pole definite real-type algebraic
curves was revealed. A discriminant curve of a subnormal operator of
finite type was defined there;
it is an algebraic curve of real type such that
all its non-degenerate irreducible pieces are real-type pole definite
(see formulas (\ref{CLa}), (\ref{discr}) below). It was
shown in \cite{Xthr}, \cite{Ysubtwo} that
a subnormal operator can be modeled as a multiplication operator by the
coordinate $z$ on a direct sum of certain (vector) Hardy classes
over the $+$ parts of the irreducible pieces of the curve. The converse
statement also is true (see \cite{Ysubtwo} and Theorem E in \S3 below).
In \S\S6, 7 of the present work, we will make use of these results. All necessary
definitions will be repeated. We also make some
small corrections to the
formulations in \cite{Ysubone}, \cite{Ysubtwo}.

\begin{thm}   
An operator $S$ is a subnormal operator of finite type without point
masses if and only if it is unitarily equivalent to a vector Toeplitz
operator $T_F$ for some symbol $F$ of class $\NDARN_m$ for some $m$.
\end{thm}

It will follow from the proof that the function $F$ is not determined uniquely.

Any scalar rational nonconstant
function $F$ on $\BD$ without poles in $\clos\BD$ belongs to
$\NDARN_1$. Since the analytic Toeplitz operator $T_F$
on $H^2$ has a normal extension, which is the operator of multiplication
by $F$ on $L^2(\BT)$, $T_F$ is subnormal. Since this normal
extension has no non-zero eigenvectors,
$T_F$ has no point masses. As it follows from formula
(\ref{Ga}), for any such $F$, $T_F$ is of finite type.
This illustrates one of implications in Theorem 2 for this simple case.

The logic of our exposition is as follows.
In \S1, we discuss a special class of algebraic curves in
$\BC^3$, which we call Ahlfors type curves.
We formulate Theorem 3, which gives a relationship
between rational matrix functions in $\NDARN_m$, pole
definite curves in $\BC^2$ and Ahlfors type curves in $\BC^3$.
In \S2, among other things,
we define vector Hardy spaces $H^2$
of a bordered Riemann surface.
In \S3, we formulate Theorem 4, which characterizes commuting pairs
of operators $(S,V)$ such that $S$ is
a subnormal of finite type without point masses and
$V$ is a pure isometry.
Theorems 1 and 4 are proven in \S4,
Theorem 2 in \S5 and Theorem 3 in \S6. The proof of Theorem 3 is
based on Theorem 4 and the existence of a pure isometry
$V$, which commutes with a given subnormal $S$ of finite type.
This is the assertion of Lemma 6, which
plays a crucial role. This lemma is derived from the structure result for
subnormal operators of finite type from
\cite{Xthr}, \cite{Ysubtwo}.

Note that Theorem 3 has purely algebraic formulation,
but our method of proving it relies on Operator Theory.

A subnormal operator of finite type is determined uniquely by two matrices.
In \S7, we calculate these matrix parameters of a subnormal operator
$S$ without point masses
in terms of a function $F$ in $\NDARNmm$ such that
$S$ is unitarily equivalent to $T_F$. In \S8, we describe a method for constructing
matrices of classes $\NDRN_m$, $\NDARN_m$. Some references to related
fields are given in the final \S9.

Some of our arguments resemble the constructions by S. Fedorov and B. Pavlov
\cite{Fed90}, \cite{Fed97}, \cite{Pav} and by Abrahamse and Bastian
\cite{AbrahBast}.

\textbf{Acknowledgements.} The author is grateful to the referee for his/her
suggestions, which helped to improve the article.

\section{Ahlfors type functions, Ahlfors type curves, and rational matrix functions
of classes \\ $\NDRN_m$, $\NDARN_m$.}

\begin{defn}
Suppose that $\De$ is a real-type algebraic curve
in $\BC^2$, all
whose irreducible pieces are nondegenerate (that is, $\wDe=\wDendeg$).
We call a function $\xi$ on $\wDe$
\textit{an Ahlfors type function} if
$\xi$ is globally meromorphic on each
irreducible piece of $\De$ and
for $\de\in\wDe$, $|\xi(\de)|=1$ if and only if
$\de$ belongs to the real part $\wDeR$ of $\wDe$.
\end{defn}

Let $\wDe_j$ be an irreducible piece of $\wDe$.
It follows that $\wDeR$ divides $\wDe_j$ into a union of
two disjoint open sets,
namely, $\{|\xi(\de)|<1\}$ and $\{|\xi(\de)|>1\}$.
Since $\wDe_j\sm\wDeR$ has at most two connected components,
these components are exactly these two subsets.
We conclude that if $\wDe$ has an Ahlfors type function,
then $\wDe$ is separated.

By the Schwartz reflection principle,
every Ahlfors type function satisfies
$\xi(\de^*)=\bar \xi(\de)^{-1}$, $\de\in\wDe$.

For any irreducible piece $\wDe_j$ of $\wDe$,
the function $\xi$, restricted
to the component $\{\de\in\wDe_j: |\xi(\de)|<1\}$
is a branched covering of the unit disc.

As proved Ahlfors in 1950, for any
compact bordered Riemann surface $\Om$ there
exists a branched covering $\xi:\Om\to\BD$, where
$\BD$ is the unit disc, see
 \cite{Ahlf}, Theorem 10.
The proof made use of a certain extremal problem.
If $\Om$ has $p$ handles and $q$ boundary contours, then \cite{Ahlf}
the degree $N$ of the extremal Alhfors function satisfies $q\le N\le 2p+q$.

Let $\wDe$ be a separated real type algebraic curve
without degenerate pieces, $\wDe_j$ one of its
irreducible pieces, and $\wDe_j^+$ (any) of the
connected components of the complement
$\wDe_j\sm\wDeR$. Let $\xi:\wDe_j^+\to\BD$ be a branched covering.
Then, by the Schwartz reflection,
$\xi$ continues to an Ahlfors type function on $\wDe$.
We come to the following

\begin{prop}
A real type algebraic curve has an
Ahlfors type function if and only if it is separated.
\end{prop}

According to our definition,
a general Ahlfors type function needs not to be
related to an extremal problem, so that its degree needs
not to satisfy the
above bounds. We refer to \cite{Fish}, Chaper 5,
and to \cite{Yam} and references therein for more information.

\begin{defn}
Let $\De_{3}$ be an algebraic curve in $\BC^3$
and $\wDe_{3}$ its desingularization.
We call it {\it an Ahlfors type curve} if

a) The coordinate functions
(denoted by $z$, $w$, $t$ it the sequel)
are non-constant on each irreducible piece of $\De_{3}$;

b) $\De_{3}$ is invariant under the
anti-analytic involution
$(z,w,t)\mapsto(\bar w,\bar z,\bar t^{-1})$;

c) Every point in $\wDe_3$ with $|t|=1$ is a fixed point of this involution.
\end{defn}

Let $\De_3$ be an Ahlfors type curve. then
the complement of the set $|t|=1$
in any of irreducible pieces of $\De_3$
is not connected. By the general theory of Klein surfaces,
each of irreducible pieces of $\wDe_3$ is separated, and
the set $|t|=1$ divides it
into exactly two components, defined, respectively,
by the inequalities $|t|<1$ and $|t|>1$.
The following properties are straightforward.

\begin{prop} 
1) If $\De_3$ is an Alhfors type curve in $\BC^3$, then
its projection onto the plane $zw$ is a real-type
separated algebraic curve in $\BC^2$. All
its irreducible pieces are non-degenerate.

2) Conversely, let $\De$ be a real-type non-degenerate
algebraic curve in $\BC^2$
without degenerate components, and let $\xi$ be an Ahlfors
type function on it. Then the graph curve
$$
\big\{
(z,w,\xi((z,w)))\in\BC^3: \quad (z,w)\in\De
\big\}
$$
is an Ahlfors type curve.
\end{prop}

If the projection of $\De_3$
onto the $zw$ plane is  pole definite, then we call
$\De_3$ a \textit{pole definite Ahlfors type curve}. We remark that
in general, the degrees of this projection on irreducible pieces of
$\De_3$ can be greater than one.
If $\De_3$ is pole definite Ahlfors type curve and
$\De_{3,j}$ is one of its pieces, then we define
its ``halves'' $\De_{3,j}^\pm$ as
in the Introduction, by requiring that
$z$ has no poles on the $+$ part. Then
$|t|<1$ on
one of the halves and
$|t|>1$ on the other.
We put $\big(\De_3\big)_\pm=\cup_j\,\De_{3,j}^\pm$.

Let $F$ be a matrix function in $\NDRNmm$ (see the Introduction).
Let $P$ be the set of poles of $F$, and put
\begin{align}
\De_{(3)}^0(F)&=
\big\{
(z,w,t)\in\BC^3:  \qquad t\ne0,\quad t,\bar t^{-1}\notin P\text{ and }\nn \\
&\qquad  \quad
\exists \phi\in \BC^m, \; \phi\ne0:
\big(F(t)-zI\big)\phi=
\big(F^*(\bar t^{-1})-wI)\big)\phi=0
\big\}.   \nn
 \nn
\end{align}
We put $\De_{(3)}(F)$ to be the closure of
$\De_{(3)}^0(F)$ in $\BC^3$ and
\beqn
\De_{(2)}(F)  =
\big\{
(z,w)\in\BC^2: \quad\exists t\in \BC: (z,w)\in \De_{(3)}(F)
\big\}
 \nn
\neqn
(so that $\De_{(2)}(F)$ is the  projection of
$\De_{(3)}(F)$ onto the $zw$-plane). It is easy to see that
$\De_{(3)}(F)$ and $\De_{(2)}(F)$ are algebraic curves.

For a fixed $t$, there are finitely many pairs $(z,w)$ such that
$(z,w,t)\in \De_{(3)}(F)$. Hence $\De_{(2)}(F)$ and
$\De_{(3)}(F)$ always have complex dimension one.

The projection $(z,w,t)\mapsto t$ maps any irreducible piece of $\De_{(3)}(F)$ onto the
whole $t$-plane $\BC$.
For all but finitely many $(z,w,t)$ in $\De_{(3)}(F)$,
$$
\Ker\big(F(t)-zI\big)=\Ker\big(F^*(\bar t^{-1})-wI\big)
$$
(because it is so for points of $\De_{(3)}(F)$
with $|t|=1$), and the dimension of these eigenspaces
is positive. This dimension defines an integer-valued multiplicity function
$\nu(\de)$ of a point $\de=(z,w,t)\in \De_{(3)}(F)$.
It follows that  $\De_{(3)}(F)$ is always an Alhfors type curve.

The desingularization of $\De_{(3)}(F)$ is a
finite union of irreducible pieces   $\wh\De_j$.
There exist positive integers $\al_j$
such that
$\nu(\de)\equiv\al_j$ on $\De_j$ (except for a finite number of points).

\begin{thm}      
Let $\De$ be an algebraic curve in $\BC^3$. Then $\De$ is a
pole definite Ahlfors type algebraic curve
such that $|t|<1$ on $\De_+$ if and only if
there is a non-degenerate analytic rational normal matrix function
$F$ such that $\De=\De_{(3)}(F)$.
\end{thm}

\begin{cornn}
An algebraic curve $\De$ in $\BC^2$ is real-type, separated and pole definite
if and only if there exist $m\ge1$ and a matrix function $F\in\NDARN_m$
such that $\De=\De_{(2)}(F)$.
\end{cornn}

\noindent This follows at once from the Theorem and from
Propositions 1 and 2.

If $F\in\NDRN_1$, then $\De_{(3)}(F)$ coincides
with the image of the map
$t\mapsto\,\big(t,F(t),\overline{F(\vvph\bar t^{-1}})\big)$,
$t\in\BC$. If, moreover, $F$ is analytic and univalent
on $\BD$, then the quadrature domain $F(\BD)$ equals to
the $z$-projection of $\De_{(3)}(F)$.

It would be interesting to see an algebraic proof of Theorem 3 and
to extend this result to functions $F$ in $\NDRN_m$. We will give an
operator theoretic proof.
The ``if'' part of Theorem 3 is straightforward.

It is common in the algebraic geometry  to consider algebraic curves
in $\BC^\ell$ as imbedded in the projective space
 $\BP^\ell(\BC)$, see \cite{FUL}. Here we do not make
an explicit use of this point of view.

\section{Pure Isometries and Vector Hardy Spaces}

We recall that an operator $V$ on a Hilbert space $H$ is called
\textit{an isometry} if $\|Vh\|=\|h\|$ for all $h\in H$.
If $V$ is an isometry, it does not follow that $V^*$ also is;
if it does, then $V$ is a unitary operator.
An isometry $V$ is called \textit{pure} if
it does not have a non-zero invariant subspace $H_1\subset H$
such that $V|H_1$ is unitary.

Let $m\ge1$.
An example of a unitary operator is given by the
operator of the forward shift operator
on the space $l^2(\BZ,\BC^m)$ of two-sided vector-valued sequences:
$$
U\{a_n\}_{n\in \BZ}=\{a_{n-1}\}_{n\in \BZ}, \qquad \{a_n\}\in\BC^m.
$$
The Fourier transform
\beqn
\F:\{a_n\}\in l^2(\BZ,\BC^m)\mapsto f(t)=\sum_{n\in\BZ}a_nt^n\in L^2(\BT,\BC^m)
\label{Four}
\neqn
is a unitary isomorphism of $l^2(\BZ,\BC^m)$ onto $L^2(\BT,\BC^m)$.
This transform is in fact a spectral representation of $U$
in the sense that
$$
\F U\F^{-1}f(t)=tf(t), \qquad f\in L^2(\BT,\BC^m).
\label{sptrl}
$$
We put $\BZ_+=\{n\in\BZ: n\ge0\}$, $\BZ_-=\{n\in\BZ: n<0\}$.
It is easy to see that
$l^2(\BZ_+,\BC^m)$
is an invariant subspace of $U$, and $V=U|l^2(\BZ_+,\BC^m)$ is a pure isometry.
In the spectral representation of $U$, operator $V$ takes the form
$$
\F U\F^{-1}f(t)=tf(t), \qquad f\in H^2_m.
\label{H2}
$$

We will also need the Hardy space
$$
H^2_{-,m}=\F l^2(\BZ_-,\BC^m).
$$
One has an orthogonal sum decomposition
$L^2(\BT,\BC^m)=H^2_{-,m}\oplus\Htwom$.
Formula (\ref{Four}) permits one to
interpret functions in $H^2_{-,m}$
as boundary values of functions $f(t)$, analytic
on $\widehat\BC\sm \clos\BD$ with $f(\infty)=0$ (here
$\wh\BC=\BC\cup\{\infty\}$  is the
Riemann sphere).

\begin{defn}
Let $H$ be a Hilbert space and $V:H\to H$ an isometry. Let $m\in\BN$.
An operator $R: H\to H^2_m$
will be called {\it a Kolmogorov-Wold representation of} $V$
if $R$ is a unitary operator, which transforms $V$ into the multiplication
operator by the independent variable:
$$
RVR^{-1}f(t) = tf(t), \qquad f\in H^2_m.
$$
\end{defn}

The following statement is a particular case of the Kolmogorov - Wold lemma, see
\cite{Nikbook}.

\begin{prop}[Kolmogorov-Wold] 
Let $H$ be a Hilbert space and $V:H\to H$ an isometry.
Then there exists a Kolmogorov-Wold representation $R: H\to H^2_m$ of
$V$ if and only if $V$ is pure and $\dim (H\ominus VH)=m$.
\end{prop}

{\bf Vector Hardy Spaces. } Let $\Om$ be a non-compact Riemann surface
with piecewice smooth boundary $\pl\Om$ (we assume that
$\clos \Om=\Om\cup\pl\Om$ is compact and is embedded into
a larger Riemann sufrace without boundary
and that $\Om$ equals to the
interior of $\clos \Om$). Then
the Dirichlet problem in $\Om$ is uniquely solvable,
that is, for any $f\in C(\pl\Om)$
there is a unique $h\in C(\clos \Om)$, which is harmonic in $\Om$ and
satisfies $h|\pl\Om=f$. Pick a point $p_0\in\Om$. Since the dual
space to $C(\pl\Om)$ equals to the space of finite Borel measures on
$\pl\Om$, it follows that there is a unique measure $d\om=d\om_{p_0}$ such that
the formula
$$
h(p_0)=\int f\,d\om
$$
holds for all functions $f$ and $h$, related as above. The measure
$d\om$ is positive and is called
{\it the harmonic measure for $\pl\Om$ at} $p_0$.

\

Assume that $p_0$ is fixed. Let $W$ be a $k \times  k$
measurable matrix valued weight on $\pl \Om$.
We say that  $W$ is \textit{admissible} if
$W>0$ and $W$, $W^{-1}$ are essentially bounded.
The corresponding weighted Hardy
class $H^2_k(W, \Om )$ consists of analytic
vector-valued functions $f:\Om\to\BC^k$ such that
the function $\|f(\cdot)\|^2$ has a harmonic majorant in $\Om$.
Each such function $f$ has boundary values a.e. on $\pl\Om$.
The norm
$$
\|f\|^2\defin\int_{\partial\Om}\langle Wf,f \rangle d\omega
$$
makes $H^2_k(W, \Om )$ a Hilbert space
(see \cite{Hasu}, \cite{Fish} or \cite{Ysubtwo}, \S9 for more details).

For any function $g\in H^\infty(\Om)$, the operator of multiplication by $g$ on
$H^2_k(W,\Om)$ is subnormal. (In general, we denote by $M_G$ the
operator of multiplication by a function $G$: $M_G f=G\cdot f$). In
particular, this applies to any bounded domain $\Om$ in $\BC$
with piecewise smooth boundary.

It is shown in the general theory of
subnormal operators \cite{Conw} that every
subnormal operator $S:H\to H$
has a \textit{minimal normal extension} $N:K\to K$,
$K\supset H$, in the sense that $N$ has no
invariant subspace $K_1$, $H\subset H_1\subsetneq K$ such that
the restriction of $N$ to $K_1$ is normal. The minimal normal extension
is unique up to the unitary equivalence.

The relationship between subnormal operators of finite type and quadrature
domains is seen from the following result.

\begin{thmC}[{\rm McCarthy, Yang \cite{McCY} }]
Let $\Om$ be a bounded finitely connected domain in
$\BC$ and $\rho$ a scalar admissible weight on $\pl\Om$. Then the operator
$$
M_zf(z)=zf(z)
$$
on $H^2(\rho,\Om)$ satisfies
$\rank (M_z^*M_z-M_zM_z^*)<\infty$ if and only if
$\Om$ is a quadrature domain.
\end{thmC}

It is also easy to show (see,
for instance, \cite{Ysubtwo}, Lemma 9.2) that the
above operator
$M_z$ on $H^2(\rho,\Om)$  is pure subnormal and its minimal normal extension
is the operator $M_z$ on $L^2(\pl\Om,d\om)$.

In the next section, we will need the following notation.
Let $\Om$ be a bordered Riemann surface as above,
let $\la\in\BC$ and let $\tau:\Om\to\BC$ be a non-constant
holomorphic function. We put
\beqn
\label{5'}
\ind_\la(\tau,\Om)=\sharp\{\de\in\Om:\; \tau(\de)=\la\},
\neqn
where the solutions $\de$ of the equation
$\tau(\de)=\la$ are counted with their multiplicities.
If $\tau $ is continuous on $\clos\Om$ and analytic in $\Om$,
then the function
$\la\mapsto\ind_\la(\tau,\Om)$ is locally constant on
$\BC\sm\tau(\pl\Om)$.

\section{Vector Toeplitz operators and subnormal operators}

Here we introduce some extra notation and
formulate Theorem 4, which will be used in the proof of
Theorem 3.

For a linear operator
$T$ on a Banach space $H$, consider an open set
$$
\rho_l(T)=
\big\{
\la\in\BC:\qquad
\exists \eps>0:\enspace \|(T-\la I)x\|\ge \eps\|x\|, \quad x\in H
\big\}.
$$
For $\la\in\rho_l(T)$, the image
$(T-\la I)H$ is closed in $H$.
The function
$\ind_T:\rho_l(T)\to\BZ_+\cup\{\infty\}$, defined by
$$
\ind_T(\la)=\dim
\big(
H\ominus (T-\la I)H
\big)
$$
is locally constant on $\rho_l(T)$.

Take any matrix-valued function $F$ in
$L^\infty(\BT,\BC^{m\times m})$. We define the
(vector) Toeplitz operator
$T_F$ in $\CL(H^2_m)$ and the
(vector) Hankel operator $\Ga_F$ in
$\CL(H^2_m,H^2_{-,m})$ by
$$
T_Fx=P_+(F\cdot x),\quad \Ga_Fx=P_-(F\cdot x).
$$
The function $F$ is called \textit{the symbol}  of these operators.
If $F\in\Hinfmm$, then $T_F=M_F$.

Let $F$ be a function in $\Hinfmm$, which is continuous on the closed unit disc
(that is, we assume that the entries of $F$ are in the disc algebra).
Choose the curves $\zeta_j$ as in Theorem 1, and let
\beqn
\label{gammaF}
\ga(F)=\bigcup_j\,\ze_j([0,2\pi])=
\bigcup_{s\in[0,2\pi]}\sigma\big(F(e^{is})\big).
\neqn
It is known \cite{BotSilb} that
for a function $F\in\Hinfmm$, continuous on the
closed unit disc,
\beqn
\label{indTF}
\rho_l(T_F)=\BC\sm\ga(F), \qquad
\ind_{T_F}(\la)=\wind_F(\la), \quad\la\in\rho_l(T_F).
\neqn
Many much more general facts about different types of spectrum
of scalar and vector Toeplitz operators are known, see, for instance,
\cite{BotSilb} and \cite{Nikbook}.

The next fact is also well-known.

\begin{prop} 
A bounded operator $T$ on $\Htwom$ commutes with the
shift operator $M_t$ on $\Htwom$ if and only if
$T=T_F$ for some symbol $F\in\Hinfmm$.
\end{prop}

It is clear that for any
normal symbol $F\in\Hinfmm$, operator $T_F$ on $\Htwom$
is subnormal (not necessarily pure), and operator
$M_F$ on $L^2_m(\BT)$  is its normal extension.

The next three lemmas
will be proven in the next section.

\begin{lem} 
A vector Toeplitz operator $T_F$
with symbol $F\in \Hinfmm$ has a finite rank self-commutator
iff $F$ is normal and rational.
\end{lem}

\begin{lem}  
Let $F\in\Hinfmm$ be a normal symbol. The operator $T_F$ is pure
subnormal if and only if
$F$ is non-degenerate. If this condition is fulfilled, then the
operator $M_F$ on $L^2_m(\BT)$ is the minimal normal extension of
$T_F$.
\end{lem}

\begin{lem}   
Let $F\in\NDARN_m$, and define
$\ga(F)$ by (\ref{gammaF}).
Let $\la\notin \gamma(F)$. Then
\beqn
\label{ind}
\wind_F(\la)=
\dim \big(\Htwom\ominus (M_F-\la)\Htwom\big)=\ind_\la\big(z,\De^+_{(3)}(F)\big),
\neqn
where the last index is defined in (\ref{5'}).
\end{lem}

\begin{lem} 
A vector Toeplitz operator $T_F$
with symbol $F\in \Hinfmm$
is a subnormal operator of finite type
without point masses
if and only if its symbol $F$ belongs to $\NDARN_m$.
\end{lem}

\begin{proof}
The ``only if'' part follows directly from Lemmas 1 and 2.
Conversely, suppose that $F\in\NDARN_m$.
Then, by the same lemmas, $T_F$
is pure subnormal and has finite rank self-commutator. It is
clear that for any $c\in\BC$, the set of points $t\in\BT$ such that
$\det(F(t)-cI)=0$ has zero measure.
Hence the operator $M_F$ on $L^2_m(\BT)$
has no non-zero eigenvectors.
By Lemma 2, this
operator is the minimal
normal extension of $T_F$. Therefore
$T_F$ has no point masses.
\end{proof}

\begin{thm}  
Let $V:H\to H$ be a pure isometry with $\codim VH=m<\infty$,
and let $S:H\to H$
be a bounded operator. Let $R: H\to \Htwom$ be
the Kolmogorov--Wold representation of $V$.
Then the following conditions are equivalent:

\textbf{(1)} $S$ is a pure subnormal operator
without point masses of finite type
and operators $S$ and $V$ commute;

\textbf{(2)} $RSR^{-1}=T_F$ for some $m\times m$ matrix symbol
$F$ of class $\NDARN_m$.
\end{thm}

Operators $S$ and $V$ that satisfy part (1) of this theorem
are a very particular
example of $n$-tuples of commuting subnormal operators.
We refer to works by Xia
\cite{Xadd}, \cite{Xtrace2} and others for a general study
of $n$-tuples of commuting subnormal and hyponormal operators.

\section{Proofs of Theorems 1 and 4 and of auxiliary lemmas}

First we will prove the theorems modulo Lemmas 1--3.
\begin{proof}[Proof of Theorem 1]
{\bf a)}
Let $\Om$ be a quadrature domain. Let $\tau=\tau(z)$ be an Ahlfors
type function in $\Om$, and fix some admissible weight
$\rho$ on
$\partial \Om$. Notice that by \cite{Ahlf}, Theorem 10,
one can find $\tau$, whose degree
equals to the connectivity number of $\Om$.
Consider operators $S=M_z$ and $V=M_\tau$ on
$H^2(\Om,\rho)$. Then $S$ is a pure subnormal operator
of finite type (see Theorem A), $V$ is an isometry, and
these operators commute.
Let $m$ be the degree of $\tau$, then it is easy to see that
$$
\dim\big(H^2(\Om,\rho)\ominus \tau H^2(\Om,\rho)\big)=m.
$$
Let $R:H^2(\Om,\rho)\to \Htwom$
be the Kolmogorov -- Wold representation of $V$. Put $T=RSR^{-1}$.
Then $T$ commutes with the shift operator
$M_t$ on $\Htwom$
(we denote by $z$ the independent variable in $\Om$
and by $t$ the independent variable in $\BD$).
 Hence $T=T_F$ for some symbol $F=F(t)\in\Hinfmm$.
It follows from Theorem C and Lemma 2 that $F\in \NDARN_m$.

Notice that $M_z-\la I$ is left invertible
for all $\la\notin \pl\Om$, and that
$$
\ind_{M_z}(\la)=
\begin{cases}
1, &\quad  \la\in \Om\\
0, &\quad \la\notin \clos\Om.
\end{cases}
$$
Since
$T_F$ is unitarily equivalent to the multiplication
operator $M_z$ on $H^2(\Om,\rho)$,
formula (\ref{indTF}) implies that assertions
{\bf 2)} and
{\bf 3)} of Theorem 1 hold. By (\ref{wind}),
{\bf 1)} is a consequence of {\bf 2)} and {\bf 3)}.

{\bf b)} Conversely, suppose that $F$ is in $\NDARNm$ and
$\Om$ is related with $F$ by means of conditions
{\bf 1)}--{\bf 3)}.
Consider the separated real-type algebraic curve $\De_2=\De_{(2)}(F)$ in
$\BC^2$. Then, by (\ref{ind}), the projection onto the $z$ plane is
one-to-one on $\De_2$ and
$\Om=z(\wDe_+)$. By Theorem A, $\Om$ is a quadrature domain.
\end{proof}

\begin{proof}[Proof of Theorem 4]
Let us prove that (1) implies (2).
Put $T=RSR^{-1}$.
Since $T$
commutes with $M_t$
on $\Htwom$, it follows that
$T=T_F$ for a matrix function
$F=F(t)\in H^\infty_{m\times m}$.
We remind that $R$ is an isometric isomorphism.
Hence $T_F$ has the same operator properties
as the operator $S$. By Lemma 4, $F\in\NDARN_m$.

The converse implication also follows easily from Lemma 4.
\end{proof}

\begin{proof}[Proof of Lemma 3]
The first equality follows from
(\ref{indTF}). Let us prove the second one.
Suppose first that $\la$ satisfies the following assumption:
the eigenspace and the root space of $F(t)$ that correspond to eigenvalue $\la$
coincide for all $t\in\BD$.
It is clear that  $\wind_F(\la)$
equals to the number of zeros of the function
$f(t)=\det\big(F(t)-\la I\big)$ in $\BD$, counted with multiplicities.
The extreme right term  equals to the sum
$$
\sum_{t\in\BD}\dim\Ker\big(F(t)-\la I\big).
$$
This implies (\ref{ind}) for points $\la$ with the above property.

For all but a finite number of points
$t$ in the unit disc, $F(t)$ has no non-trivial
Jordan blocks. Hence only a finite number of values of $\la$
were excluded from our consideration.
Since all terms in (\ref{ind}) are locally constant on $\BC\sm\ga(F)$,
the general case follows.
\end{proof}

\begin{proof}[Proof of Lemma 1]
We will use the known formula
\beqn
\label{Ga}
T_F^*T_F-T_FT_F^*=\Ga_{F^*}^*\Ga_{F^*}+T_{F^*F-FF^*}, \qquad F\in\Hinfmm.
\neqn
To prove it, notice that
for any $G\in \Linfmm$, $T_G^*=T_{G^*}$, and
$$
T_{G^*G}-T_G^*T_G=\Ga_G^*\Ga_G.
$$
Then (\ref{Ga}) is obtained by putting
$G=F$ and $G=F^*$ and taking into account that $\Ga_F=0$.

Fix any $x=x(t)$ in $H^2_m$. It is easy to see that $t^nx$ tends
weakly to zero and that $\|\Ga_{F^*}(t^nx)\|\to0$ as $n\to\infty$. By
general properties of Toeplitz operators (see \cite{Nikbook})
and (\ref{Ga}), $\|(F^*F-FF^*)x\|=\lim_{n\to\infty}\|T_{F^*F-FF^*}t^nx\|=0$. Therefore
$F^*F-FF^*\equiv0$ a.e. on $\BT$, so that $F$ is a normal symbol.
Applying (\ref{Ga}) once again, we get that $\Ga_F$ is a finite
rank operator. By Kronecker's lemma, $F$ is rational
(see \cite{Nikbook}, p. 183 for the scalar case; the vector case
has the same proof).
\end{proof}

Before proving Lemma 2, we will need one more fact.

\begin{lem} 
Let $F$, $\theta$, $G$ be $m\times m$, $m\times s$ and $s\times s$
constant matrices,
respectively, where $1\lews s\lews m$. If $F$ is normal,
$\theta^*\theta=I$ and $F^*\theta=\theta G$, then $G$ is normal and $F\theta=\theta G^*$.
\end{lem}

\begin{proof}
For any polynomial $p$, $p(F^*)\theta=\theta p(G)$. Since $F$ is normal,
we can find a  polynomial $p$ with simple roots such that $p(F^*)=0$. It
follows that $G$ also has  no non-trivial Jordan blocks.

Let $a_1,\dots,a_s\in\BC^s$ be a complete family of eigenvectors of $G$
and $\la_1,\dots,\la_s$ the corresponding eigenvalues. We assume that
if $\la_j=\la_k$ and $j\ne k$, then $\langle a_j,a_k\rangle=0$.

We have
$$
(F^*-\la_j I)\theta a_j=\theta(G-\la_j I)a_j=0.
$$
Therefore
if $\la_j\ne\la_k$, then
$\langle a_j,a_k\rangle=\langle \theta a_j,\theta a_k\rangle=0$. Hence
$a_1,\dots,a_s$ form an orthogonal basis, so that $G$ is normal. From
the normality of $F$ and $G$ we get
$$
F\theta a_j=\bar\la_j\theta a_j=\theta G^*a_j,\quad j=1,\dots,s,
$$
which implies that $F\theta=\theta G^*$.
\end{proof}

Recall that a function $\Phi(t)\in \Hinfms$ is called \textit{inner}
if its boundary values on $\BT$ are isometries a.e.
(then it follows that $s\le m$).

\begin{proof}[Proof of Lemma 2]
Suppose first that
$F$ is degenerate, that is,
there is a constant $c\in\BC$ such that
$\det(F(t)-cI)\equiv0$ for $t\in\BD$. Then
the same holds for a.e. $t\in\BT$. Put
\beqnay
L&\defin\bigl\{x\in\Htwom:Fx=cx\text{ a.e. on }\BT\bigr\}  \nn  \\
&\,=\bigl\{x\in\Htwom:F^*x=\bar cx\text{ a.e. on }\BT\bigr\}.  \nn
\neqnay
Note that $(F-cI)x\equiv0$ a.e. on $\BT$ iff
$(F-cI)x\equiv0$ in $\BD$. Put $X=(F-cI)^\sim$, where $(F-cI)^\sim$ is the
transpose associate matrix to $F-cI$. Then $X\in\Hinfmm$ and
$(F-cI)X\equiv0$.
If $X\not\equiv0$, then $L\ne0$.
If $X\equiv0$, then there exists an integer $k$, $0\lews k<m$
and a $k\times m$ submatrix $G$ of $F-cI$ such that for a vector $x$ in
$\Htwom$, $(F-cI)x\equiv0$ in $\BD$ if and only if $Gx\equiv0$ in $\BD$.
We take the smallest possible $k$.
Let $G'$ be any matrix in $H^\infty_{(m-k)\times m}$ such $\det J\not\equiv0$,
where $J\defin\left(\smma G\\G'\esmma\right)$. Then it is easy to see
that $x\defin J^\sim(0,\dots,0,1)^{\text{\smc T}}\not\equiv0$ is in
$L$.

It follows that in all cases, the subspace $L$ is non-zero, closed and
is invariant both for $T_F$ and $T_F^*=P_+M_{F^*}$. Hence $T_F$ is not
pure.

Let us prove the converse.
Suppose that $F$ is non-degenerate.
In general, if $S\in\CL(H)$ is a subnormal
operator (not necessarily pure) and $N\in\CL(K)$ its normal extension,
then the maximal
invariant subspace of $S$ on which $S$ is normal is given by
$
H_1=\bigl\{x\in H: N^{*k}x\in H\quad\forall k\in\BN\bigr\}
$
(see \cite{Conw}). So in our case,
$$
H_1=\bigl\{x\in \Htwom: F^{*k}x\in \Htwom\quad\forall k\in\BN\bigr\}
$$
and we have to prove that $H_1=0$. We remark that $H_1$ is a closed
$M_z$-invariant subspace of $\Htwom$. If $H_1\ne0$, then by the
Beurling--Lax--Halmos theorem (see \cite{Nikbook}), there exists a
natural $s$, $1\lews s\lews m$ and a matrix function
$\theta\in H^\infty_{m\times s}$
such that $\theta(t)$ is an isometry for a.e. $t\in\BT$ and
$H_1=\theta\Htwos$. For all $r\in\BC^s$, $F^*\theta r\in\Htwos$. Therefore there
exists a function $G$ in $\Hinfss$ such that
\beqn
\label{t13.1}
F^*\theta=\theta G\quad \text{ a.e. on }\BT.
\neqn

For any complex $c$, we can replace $F$, $G$ by $F-cI$, $G-\bar cI$. We
will assume without loss of  generality that $\det G(0)=0$.

Lemma 5 yields
\beqn
\label{t13.2}
F\theta=\theta G^*\quad \text{ a.e. on }\BT.
\neqn
The matrix $\theta$ has a $s\times s$ minor whose determinant is not
identically zero on $\BT$. Therefore there exists a constant
$s\times m$ matrix $\rho$ such that $\det(\rho\theta)\not\equiv0$.
Note that $\det(\rho\theta)\in \Hinf$.  By
(\ref{t13.2}),
\beqn
\label{t13.3}
\rho F^n\theta=\rho\theta G^{*n}\quad \text{ a.e. on }\BT
\neqn
for all $n\in\BN$. Put $g=\det G$, then $g\in H^\infty$. By
(\ref{t13.3}),
\beqn
\label{t13.4}
\bar g^n\det(\rho\theta)|\BT\in H^\infty
\neqn
for all $n\in\BN$. We obtain from the Nevanlinna factorization
that $\bar g=\phi^{-1}h$
on the unit circle, where $h\in H^\infty$ and $\phi$ is inner.
It follows that for every $n$,
the function $h^n\det(\rho\theta)$ has an inner multiple
$\phi^n$. This implies that $\phi$ divides $h$.
Therefore $\bar g=h_1$
on $\BT$ for some $h_1$ in $H^\infty$, so that $g=\const$.
Since $g(0)=0$,
we get that $g\equiv0$.

Let $a(t)\in\Ker G(t)$ for $t\in\BT$, $a(t)\ne0$, a.e.
$t\in\BT$. By (\ref{t13.1}), $F^*(\theta a)\equiv0$ on $\BT$,
which yields $\det F(t)\equiv0$. This contradicts to our assumption that
$F$ is non-degenerate. We have proved that $H_1=0$, that is, that
$S$ is pure.

At last, suppose that $\det(F-cI)\not\equiv0$ on $\BT$ for all
$c\in\BC$, and let us check that $M_F$ is a \textit{minimal} normal extension of
$T_F$. We have to prove that the subspace
$$
K\defin\spa\{F^{*n}\Htwom:n\gews0\}
$$
of $L^2_m(\BT)$ coincides with $L^2_m(\BT)$. Suppose that for
some $y$ in $L^2_m(\BT)$, we have $\langle y,F^{*n}x\rangle=0$ for
all $x$ in $\Htwom$ and all $n\gews0$. Then $y$ is in $\Htwomin$.

The formula $\wh x(t)=\bar t x(\bar t)$ defines a symmetry on
$L^2_m(\BT)$, which maps  $H^2_m$
onto $H^2_{-,m}$ and $H^2_{-,m}$ onto $H^2_m$.
Put $\widehat F(t)=F^*(\bar t)$, then
$\widehat F$ is also in $H^\infty_{m\times m}$.
Applying the symmetry $x\mapsto \hat x$, we get
$\langle \widehat F^{*n}\hat y,u\rangle=0$ for all $u\defin \hat x$ in
$\Htwomin$. Therefore $\widehat F^{*n}\hat y\in\Htwom$ for all
$n\gews0$. Since $\widehat F$ has the same properties as $F$, we
conclude from the above that $\hat y=0$. Hence $K=L^2_m(\BT)$.
\end{proof}




\section{Discriminant curve of a
subnormal operator. Proof of  Theorem 2} 
Let $S$ be a subnormal operator of finite type.
Put
\beqn
\label{CLa}
\begin{aligned}
M&\defin\Range S^*S-SS^*;  \\
C&=C(S)\defin S^*S-SS^*|M, \qquad \La=\La(S)=\bigl(S^*|M\bigr)^*.
\end{aligned}
\neqn
It is known that $C>0$ and $S^*M\subset M$.
In  \cite{Xone}, \cite{Xtwo},
Xia discovered the role of operators $C$, $\La$ in the study of the
spectral structure of the operator $S$.
In \cite{Xone}--\cite{Xfive},
he constructed and studied an analytic model of a subnormal operator with
the help of these operators and a certain projection-valued function, analytic
outside the spectrum of the minimal normal extension of $S$
(``Xia's mosaic'').
One of the consequences of Xia's results is that
the pair $(C,\La)$ of operators on $M$
completely determines a pure subnormal operator $S$.

In our context of
the study of
subnormal operators of finite type
we associate with any operators $C=C^*$ and  $\La$ on
a finite dimensional space $M$
\textit{the discriminant surface}, given by
\beqn
\label{discr}
\Del=\bigl\{(z,w)\in \BC^2:\det\bigl(C-(w-\La^*)(z-\La)\bigr)=0\bigr\}
\neqn
It always is an algebraic curve of real type.

If in (\ref{discr}),
$C$ and $\La$ correspond to a finite type
subnormal operator $S$, then we will write
$\De=\De(S)$.

In \cite{Ysubone}, conditions on $\Del$
that are necessary and sufficient for the existence of
$S$ with $C=C(S)$, $\La=\La(S)$ were given. The formulations
in \cite{Ysubone} contain certain inaccuracies.
The corrections are as follows.

Let $\De$ be given by (\ref{discr}). Define, as in \cite{Ysubone},
a meromorphic $\CL(M)$-valued function $Q$ by
$$
Q(\de)\defin\Pi_w(C(z-\La)^{-1}+\La^*),
\qquad \de=(z,w)\in\De\setminus z^{-1}\bigl(\sigma(\La)\bigr),
$$
where $\Pi_w(A)$ is the Riesz projection onto the
root space of a matrix $A$ corresponding to the eigenvalue $w$.
The values $Q(\de)$, $\de\in\De$ are
parallel projections in $M$.
Let $\De_s$ be the (finite) set of singularities of $\De$.
Then
Theorem 1 in \cite{Ysubone} has to be formulated
as follows (we conserve the numeration of formulas of
\cite{Ysubone}).

\begin{thmD}[\cite{Ysubone}, Theorem 1]
Let $M$ be a finite-dimensional Hilbert space
and $C,\,\La$ operators on $M$
with $C>0$. Define $\De, Q$ as above.
Then there exists  a
subnormal operator $S$ satisfying $C=C(S)$ and $\La=\La(S)$ if and only if
the following conditions hold:

i) $\De$ is separated and pole definite;

ii) Put
\beqn
\label{4.1}
\mu(z)=\sum_{w:(z,w)\in\De_+}Q\bigl((z,w)\bigr),\quad
    z\in\mathbb C\setminus\bigl(\sLa\cup\ga\cup z(\De_s)\bigr).
\tag{4.1}
\neqn
Then there exists a positive
${\cal L}(M)$-valued measure $de(\cdot)$ such that
\beqn
\label{4.2}
(\La-z)^{-1}(1-\mu(z))=\int\frac{de(u)}{u-z},\quad
       z\in\BC\setminus\bigl(\sLa\cup\ga\cup z(\De_s)\bigr)
\tag{4.2}
\neqn
and
\beqn
\label{4.3}
\bigl(C-(\bar u-\La^*)(u-\La)\bigr)de(u)\equiv0.      \tag{4.3}
\neqn
If (i), (ii) hold, then the measure $de(\cdot)$ is connected with the
operator $S$ by the formula (1.1) (from \cite{Ysubone}), and
$\mu$ is Xia's mosaic of $S$.
\end{thmD}

In Theorem 2 in \cite{Ysubone}, one just has to replace
item i') by the following:

\

{\it i') $\De$ is separated and pole definite. }

\

Proposition 1 in \cite{Ysubone} is erroneous, namely, it may happen
that $\De$ is separated and pole definite, but the set
$\{\big|\frac{dw}{dz}\big|=1\}$ is strictly larger than
$\wDeR$. One has to define $\wDe_+$ and $\wDe_-$ only for separated
pole definite curves $\De$ (as we do in the present paper). Then
in the proof of Theorem D formula (5.10) (see \cite{Ysubone})
can just be taken as a definition of $\wDe_\pm$. It follows from
\cite{Ysubone}, Lemma 1 that $\De$ is pole definite.

We remark that for a subnormal operator $S$ of finite type,
$\De(S)$ is always separated, but can have degenerate pieces (even if
$S$ has no point masses).

The paper \cite{Ysubtwo} explains how to
construct the operator $S$, starting from
the corresponding $C$ and $\La$; this paper
was based on previous work by Xia.
The model of the operator $S$ was formulated
in \cite{Ysubtwo} in terms of weighted
analytic functional classes $H^2$ of the $+$ halves of the
components $\wDe_j$ of the curve
$\De(S)$. Slight modifications are to be done also in this paper.
Namely, in Lemma 11.4 and
Theorem 11.5 one has to replace the word ``separated'' by
``separated pole definite''. Then
Lemma 11.1 is not necessary.

Let $N$ be the minimal normal extension of $S$.
As Xia proves, the spectrum of the minimal normal extension $N$ of $S$
coincides with the $z$-projection of the real part $\wDeR$
of the discriminant curve (\ref{discr}) (with a possible exception of a
finite number of points). That is,
\beqn
\label{siN}
\si(N)\cong \{z\in \BC:
\det\bigl(C-(\bar z-\La^*)(z-\La)\bigr)=0\bigr\};
\neqn
here $A\cong B$ means that sets $A$, $B$ differ in a finite
number of points.

In what follows, we repeat briefly the results of \cite{Xthr}, \cite{Ysubtwo} that
will be used in the sequel.

Let $\De=\De(S)$,
then $\De$ is a separated pole definite
real type algebraic curve.
Let $\Dendeg=\bigcup \De_j^\kj$
($\De$ may have degenerate components).
Suppose we have admissible matrix-valued weights $W_j$
on the boundaries $\pl\wDe_j^+$.
For simplicity of notation, we denote by $W$ the
collection of matrix weights $(W_1,\dots,W_k)$, and
put
$$
H^2(W,\wDendegp)=\oplus_j\,
H^2_{\kj}(W_j,\wDe_j^+).
$$
For any choice of a non-degenerate
separated pole definite curve $\Dendeg$ and a weight $W$,
the multiplication operator
$$
(M_zf)(\de)=z(\de)f(\de)
$$
is pure subnormal of finite type; moreover, its discriminant surface
coincides with $\bigcup_j \De_j^\kj$ (\cite{Ysubtwo}, Lemma 11.4).
Subnormal operators that are unitarily equivalent to
these ones were called \textit{simple} in \cite{Ysubtwo}.

\begin{thmE}[see \cite{Xthr}, \cite{Ysubtwo}]
Let $S$ be a subnormal operator without point masses.
and let
$ \bigcup\De_j^\kj$ be the non-degenerate
part of its discriminant surface $\De(S)$.
Then there are $\kj\times \kj$
admissible matrix weights $W_j$ on
$\partial\wDe_j^+$ and a subspace
$$
\wH_1\subset H^2(W,\wDendegp)
$$
of finite codimension such that $\wH_1$ is invariant
under operator $M_z$ and $S$ is unitarily equivalent to
operator $M_z$, restricted to $\wH_1$.

Conversely, for any separated pole definite
curve $\De=\prod \De_j^\kj$ without degenerate components, any
$\kj\times \kj$ matrix weights $W_j$ and any subspace
$\wH_1$ of $H^2(W,\wDendegp)$ with the above properties,
the operator $M_z$ on $\wH_1$ will be subnormal of finite type, and the
non-degenerate part of the discriminant surface $\De(M_z)$
will be exactly equal to $\bigcup \De_j^\kj$.
\end{thmE}

A subspace $\wH_1$ of $H^2(W,\wDendegp)$ has the above two properties
iff it has a form
\beqn
\label{H1}
\wH_1=\bigl\{x\in H^2(W,\wDendegp):\quad
\langle x,\psi^j_{\la_k}\rangle=0,
\quad1\lews k\lews r,\,0\lews j\lews m_k
\bigr\}
\neqn
where (not necessarily distinct) points
$\la_k$, $1\le k\le r$ belong to $\bigcup_j z(\De_j^+)$ and
$\{\psi^j_{\la_k}\}_{j=0}^{m_k}$ are
corresponding Jordan chains of generalized eigenvectors:
$\bigl(M_z^*-\bar\la_k\bigr)\psi^0_{\la_k}=0$,
$\bigl(M_z^*-\bar\la_k\bigr)\psi^j_{\la_k}=\psi^{j-1}_{\la_k}$,
$j=1,\dots,m_k$. See \cite{Ysubtwo}, Theorem 12.3.

\begin{lem}   
For every subnormal operator $S$ of finite type without
point masses, there exists an isometry $V$ as in Theorem 4
such that $SV=VS$.
\end{lem}
 \begin{proof}
We apply Theorem E.
Let $\De$ be the discriminant surface of $S$,
and fix an Ahlfors type function $\phi$ on
its non-degenerate part $\wDendeg$ such that
$|\phi|<1$ on $\wDendegp$.
Replacing $S$ by a unitarily equivalent operator, we can assume that
$Sf(\de)=z(\de)f(\de)$, $f\in \wH_1$, where $\wH_1$ is a subspace of
$H^2(W,\De_+)$ of finite codimension for some admissible matrix weights
 $W_j$. Representation (\ref{H1})
implies that there exists a natural $N$ such
that
\beqn
\label{t13.5}
\Psi^N\cdot H^2(W,\wDendegp)\subset \wH_1\subset  H^2(W,\wDendegp),
\neqn
where
$$
\Psi(\de)=\prod_k\bigl(z(\de)-\la_k\bigr).
$$
Consider the finite Blaschke product
$B(\xi)=\prod_j\frac{\xi-\phi(\de_j)}{1-\overline \phi(\de_j)\xi}$,
where $\{\de_j\}$ are all points of $\De_+$ whose $z$-projection
coincides with one of $\la_k$. Set $\phi_1=B^{N_1}\circ\phi$, where
$N_1$ is a natural number. If $B$ is constant,
then we put $\phi_1=\phi$.
Then $\phi_1$
is also an Ahlfors type function on $\wDendeg$.
If $N_1$ is large enough, then, moreover,
(\ref{t13.5}) implies that
$$
\phi_1\wH_1\subset\Psi^N  H^2(W,\wDendegp) \subset \wH_1.
$$
In particular, the operator of multiplication by
$\phi_1$ acts on $\wH_1$. Denote this operator by $V$.
It is obviously a pure isometry, and the codimension
of $V\wH_1$ in $\wH_1$ is finite (every function in $\wH_1$,
which has zeros of sufficiently high order in zeros
of $\phi_1$ belongs to $V\wH_1$). The equality
$SV=VS$ holds, because both are multiplication operators
by scalar functions.
\end{proof}

\begin{rem} If $S$ is simple in the sense of \cite{Ysubtwo},
then one can take $\wH_1=H^2(W,\wDendeg)$. It follows that
in this case one can put $\phi_1=\phi$. The degree of $\phi_1$
on each piece $\wDe_j$ of $\wDendeg$ does not exceed $2p_j+q_j$, where
$p_j$ stands for the number of handles and $q_j$
stands for the number of boundary contours of
$\wDe_j$. If $S$ is not simple, then the minimal
possible degrees of
$\phi_1$ on irreducible pieces of $\wDendeg$ can be much higher.
\end{rem}

\begin{proof}[Proof of Theorem 2]
(1) If $F\in\NDARN_m$, then by Lemmas 1 and 2, $T_F$ is a subnormal operator
of finite type without point masses, and the same is true for
any operator unitarily equivalent to $T_F$.

(2) Conversely, let $S$ be a subnormal operator of finite type
without point masses. By Lemma 6, there exists a pure isometry $V$ that
commutes with $S$. Now Theorem 4 provides a desired
matrix symbol $F$ in $\NDARN_m$ (for some $m$) such that
$S$ and $T_F$ are unitarily equivalent.
\end{proof}

The construction of the isometry $V$ is far from unique.
Hence the symbol $F$ in Theorem 2 is also determined in a non-unique way.

\section{Proof of Theorem 3}
Let $\Pzt$ be the projection of $\BC^3$ onto its
coordinate subspace $zt$: \linebreak 
$\Pzt(z,w,t)=(z,t)$,
$(z,w,t)\in\BC^3$.

\begin{lem}  
Let $\De_3$ be an Ahlfors type curve in $\BC^3$ (not necessarily
pole definite). Then there is finite subset $\Phi$ of $\De_3$ such that
$\Pzt$ is one-to-one on $\De_3\sm \Phi$.
\end{lem}

\begin{proof}
The image $\Pzt\De_3$ is an algebraic curve in $\BC^2$.
Let $\wh{\Pzt\De_3}$ and $\wDe_3$ be the desingularizations
of the curves $\Pzt\De_3$, $\De_3$. Then $\Pzt|\wDe_3$ is a
branched covering of $\wh{\Pzt\De_3}$. The number of preimages
of a point of $\wh \De_3$ under this covering
(counted with multiplicities) is constant on each irreducible
piece of $\wh {\Pzt\De_3}$. Take an irreducible piece $K$
of $\wh {\Pzt\De_3}$. It is a projection of (at least one)
irreducible piece of $\wDe_3$. Hence the
set of solutions of the equation
$|t|=1$ on $K$ is a finite union of closed curves, in particular,
it is infinite. Any generic point
$(z,t)$ of this set has only one preimage on
$\wDe_3$, namely, $(z, \bar z, t)$
(see the definition of an Ahlfors type curve in \S1).
Therefore, in general, all but finite number of points
of $K$ have only one preimage on
$\wDe_3$. The assertion of Lemma
follows.
\end{proof}

It follows from this lemma that every
Ahlfors type curve $\De_3$
restores in a unique way from its projection
onto the  plane $zt$.

\begin{lem}   
Let $r(\cdot,\cdot)$ be a polynomial in two variables
and $F$ be a matrix function in $\NDRN_m$ for some $m$. Then
$r(M_F,M_t)=0$ if and only if $r(z,t)\equiv0$ on $\De_{(3)}(F)$.
\end{lem}

\begin{proof}
Put $\De_3=\De_{(3)}(F)$.
For $t\in\BT$,
\beqn
\label{ker}
\BC^m=\bigoplus_{(z,t)\in\Pzt\De_3}
\Ker\big(F(t)-zI\big).
\neqn
It follows that the same is true for all but a finite number of points
$t\in\clos\BD$. Consider the vector bundle $\F$
over the open subset
$\{|t|<1\}$ of $P_{zt}\De_3$ with fibers
$$
\F((z,t))=\Ker\big(F(t)-zI\big), \qquad |t|<1
$$
(the dimension $k_j$ of the fiber can be different on different
irreducible pieces of $\Pzt\De_3$).

For any meromorphic cross-section $\eta$ of $\F$, put
$$
\big(P\eta\big)(t)=\sum_{z:\;\; (z,\,t)\in\Pzt\De_3}
\; \eta\big((z,t)\big), \quad |t|<1.
$$
It follows from (\ref{ker}) that for any $h\in\Htwom$, there exists
a unique meromorphic cross-section $h^\sharp$ of $\F$ such that
$$
Ph^\sharp=h.
$$
The function $h^\sharp$ can have poles in points
$(z,t)$ such that $|t|<1$ and (\ref{ker}) is violated in $t$; the orders
of these poles are bounded by a constant
that only depends on the geometry of $\F$.  If
$h\not\equiv0$, then $h^\sharp\not\equiv0$.

It is easy to see that
\beqn
\label{split}
r(M_F,M_t)h=r(M_F,M_t)Ph^\sharp=
P\big(r(z,t)h^\sharp\big), \qquad h\in\Htwom.
\neqn
Therefore $r(z,t)\equiv0$ on $\De_{(3)}(F)$ implies that
$r(M_F,M_t)\equiv0$. Conversely, if
$r(z,t)\not\equiv0$, then (\ref{split}) implies that
$r(M_F,M_t)h\not\equiv0$ for any non-zero $h$ in
$\Htwom$.
\end{proof}

\begin{proof}[Proof of Theorem 3]
(1) Let $m\ge1$ and $F\in\NDARNm$. Then
$\De_{(3)}(F)$ is an Ahlfors type curve.
Eigenvalues of matrices $F(t)$, $|t|\le1$ are uniformly bounded.
Therefore $\De_{(3)}(F)$ is a pole definite
Ahlfors type curve and $z$ is bounded on the
subset $\{|t|<1\}$ of $\De_{(3)}(F)$. It follows that the $+$ part of
$\De_{(3)}(F)$ coincides with the subset of
$\De_{(3)}(F)$ where $|t|<1$.

(2) Conversely, let $\De$ be an Ahlfors type curve in $\BC^3$,
meeting the hypotheses of the theorem. Decompose it into irreducible curves:
$\De=\bigcup\De_j^{k_j}$.
We can define ``halves''
$\wDe^j_+$ of irreducible pieces
$\wDe_j$ of $\wDe$ so that on
$\wDe^j_+$, $z$ is bounded and $|t|<1$.
Consider the functional class
$$
H^2(W, \wDe_+)=\bigoplus_{j=1}^N H^2_{k_j}(W_j,\wDe^j_+),
$$
where $k_j\times k_j$
admissible matrix weights $W_j$ are chosen in an arbitrary way.
Consider (bounded) multiplication operators
$S=M_z$ and $V=M_t$ on $H^2(W, \wDe_+)$. Then $V$ is a pure isometry and $S$ is a pure
subnormal operator and $VS=SV$. By Theorem~4, there exists an integer $m\ge1$,
a matrix function $F\in \NDARNmm$ and
an isometric isomorphism
$R:H^2(W, \wDe_+)\to \Htwom$
such that $RSR^{-1}=M_F$, $RVR^{-1}=M_t$.
We are going to prove that $\De_{(3)}(F)=\De$,
that is, that these curves have the same irreducible pieces,
and that their multiplicities also coincide.

Let us use the notation of
Lemmas 7--8. Take any polynomial $r(z,t)$ in two variables.
Applying Lemma 8 and the above isomorphism,
we get that $r$ vanishes on $\De$ if and only if $r(M_F,M_t)=0$
if and only if $r$ vanishes on $\De_{(3)}(F)$. This implies that
$\De_{(3)}(F)$ consists of the same irreducible components as $\De$:
$\De_{(3)}(F)=\De_j^{k'_j}$ for some numbers $k'_j\ge 1$.

Choose a complex constant $\al\in\BC$ such that
the images of the real parts of the irreducible pieces of $\De$
under the function $z+\al t$ all are different. It is easy to show that
it is possible (for any two fixed pieces, the set of $\al$'s such that
these images coincide has empty interior). We can also suppose that $\al$
is chosen so that the matrix function $F(t)+\al tI$ is non-degenerate;
then $F(t)+\al tI\in\NDARNm$.

We apply Lemma 3 to this matrix function. Since
$R$ transforms the pair of operators $(M_z,M_t)$ on
$H^2(W, \wDe_+)$ into the pair $(M_F, M_t)$ on $H^2_m$, we get
$$
\ind_\la(z+\al t, \wDe_+)
= \codim \big( (M_{F(t)+\al tI})\Htwom \big)
=\ind_\la \big(z+\al t, \wDe_{(3)}^+(F)\big)
$$
for all $\la \notin (z+\al t) \big(\pl \wDe_+\big)=
(z+\al t) \big(\pl \wDe_{(3)}^+(F)\big)$.
By comparing the jump of these indices on
$(z+\al t) $- images of the components of the curve
$\pl \wDe_{(3)}^+(F)$, we deduce that $k_j=k'_j$, $j=1,\dots,N$.
\end{proof}

\begin{rems}
{\bf(1)} Suppose that
operators $S$ and $T_F$ as in Theorem 2 are unitarily equivalent.
It can be proved by the same argument as above that
$$
\De_2(F)=\Dendeg(T_F)=\Dendeg(S),
$$
and multiplicities of irreducible components are equal.

{\bf(2)}
Recall that if $S:H\to H$ is a
subnormal operator and $N$ its minimal normal extension,
then $S'=N^*|K\ominus H$ also is subnormal; this operator is called
\textit{dual to} $S$ \cite{Conw}.

Let $S:H\to H$ be a subnormal operator of finite type
without point masses and $N:K\to K$ its minimal normal extension. Let
$T_F:H^2_m\to H^2_m$ be a Toeplitz operator unitarily equivalent to
$S$ as in Theorems 2 and 4 and $R:H\to H^2_m$ the corresponding isometric
isomorphism such that $RSR^{-1}=T_F$. Then $R$ extends to a unitary
isomorphism $\wtU:K\to L^2_m(\BT)$ such that
$\wtU N\wtU^{-1}=M_F$.
One has $\wtU H'=\Htwomin$ and  $\wtU N^*=M_{F^*}\wtU$.
It follows that the dual
operator $S'$ is unitarily equivalent to $T_{\widehat F}$,
where $\widehat F(t)=F^*(\bar t)$.
\end{rems}

It follows from the remark after the proof of Lemma 6
and from the above proof that if $\De$ is
a pole definite Ahlfors type curve in $\BC^3$, then
$\De=\De_{(3)}(F)$ for some $F$
in $\NDARN_m$, with
$$
m\le \sum_j k_j(2p_j+q_j),
$$
where $\wDe=\bigcup \wDe_j^{k_j}$ and
$p_j$, $q_j$ denote the number of handles and of
boundary contours of $\wDe_j^+$, respectively.

\section{Characterization of matrix parameters}

Here we
will describe matrix parameters  $(C,\La)$
of a finite type subnormal operator $S$
\textit{without point masses}
(see (\ref{CLa})) in terms of
a matrix symbol $F$ such that $T_F$
is unitarily equivament to $S$.
First let us discuss Blaschke- Potapov products.

A \textit{Blaschke factor} is a scalar function
of the form
$b(t)=\xi\,\frac{t-a}{1-\bar a t}$, where $a\in\BD$ and
$\xi\in\BT$ are constants.

A matrix  function $B\in\Hinfmk$ is called {\it inner}
if $B(t)$ is an isometry for a.e. $t\in \BT$ (then it follows that
$k\le m$). It is known \cite{Potapov}
that an $m \times m$  matrix function $B$ is
rational and inner iff it
can be represented as
$$
B(t) = v\prod_{n=1}^M \big({b_n(t)P_n + (I-P_n)}\big),
$$
where $v$ is an $m \times m$ unitary constant matrix, $b_n$ are
Blaschke factors, and $P_n$ are orthogonal projections in $\BC^m$.
Matrix functions $B$ of this class are
called {\it finite Blaschke--Potapov products.}
In particular, a scalar function $B$
is rational and inner iff it is a finite Blaschke product:
$
B(t) = v\prod_{n=1}^M b_n(t),
$
where $v$ is a complex unimodular constant.

\begin{defn}
Let $\al$, $h$ be rational
matrix functions in $\Hinfmm$. We call these functions
\textit{right coprime} if
equalities
$\al=\al_1B$, $h=hB$, where
$B$ is a finite Blaschke--Potapov product
and $\al_1$, $h_1$ are in $\Hinfmm$
imply  that $B$ is a unitary constant.
\end{defn}

Assume that $\det\al\not\equiv0$,
$\det h\not\equiv0$ in $\BD$. It is easy to see that
in this case $\al$, $h$ are right coprime iff
$$
\Ker \al(t)\cap \Ker h(t)=0, \qquad t\in\BD.
$$

\begin{lem} 
\label{LEM9}
Let $G$ be a rational $m\times m$ matrix function such that
$\det G\not\equiv 0$.

1) There is an $m \times m$ finite Blaschke--Potapov product
$\al$ such that $\Ker \Ga_G=\al \Htwom$.

2) The above representation of
$\Ker \Ga_G$ is equivalent to a factorization
$$
G=h\al^{-1},
$$
where
$h,\al\in \Hinfmm$ are rational, $\al$
is a finite Blaschke--Potapov product and $h,\al$ are right coprime.

3) The factorization
$G=h\al^{-1}$ of the above form is unique, up to a
substitution
$h\mapsto hu$,
$\al\mapsto \al u$, where $u$ is a unitary constant.
\end{lem}

We remark that the factorization
$G(t)=h(t)\al^{-1}(t)$ in $\BC$ is equivalent to
$G(t)=h(t)\al^*(t)$, $t\in\BT$.

\begin{proof}
{\bf 1)}
It is a standard fact that $\Ker \Ga_G$ is invariant under the
shift operator $x=x(t)\mapsto tx(t)$ on $\Htwom$. By the
Beurling-Lax-Halmos theorem \cite{Nikbook}, there is an
integer $k$, $0\le k\le m$
and a matrix inner function $\al$
of size $m\times k$ such that
$\Ker \Ga_G=\al \Htwok$. In our case
of rational $G$, it is easy to find a
finite scalar Blaschke product $\phi$ such that
$\Ker \Ga_G\supset\phi \Htwom$.
It follows that $\phi I=\al\be$, where $\be$ is a
matrix inner function of size $k\times m$. Therefore $k=m$,
$\Ker \Ga_G$ has finite codimension in $\Htwom$, and $\al$ is a
finite Blaschke--Potapov product.

{\bf 2)} Let $\Ker \Ga_G=\al \Htwom$. Then
$G\al\Htwom\subset\Htwom$, hence $h\defin G\al\in\Hinfmm$. If $h$ and
$\al$ were not right comprime, that is,
$h=h_1B$,
$\al=\al_1B$ for a nonconstant rational inner function
$B\in \Hinfmm$, then $G=h\al^{-1}=h_1\al_1^{-1}$ would give
$\Ker \Ga_G\subset\al_1 \Htwom$, a contradiction.

Conversely, the same arguments show that a right coprime factorization
$G=h\al^{-1}$ implies that $\Ker \Ga_G=\al \Htwom$.

Statement 3) follows from the Beurling-Lax-Halmos theorem.
\end{proof}

For any $m$ and any matrix function $F$ in $\NDARNm$,
\enspace
$\det F$ does not vanish identically.
Note that $F^*(t)$ coincides for $t\in\BT$ with a rational
matrix function, namely, with the function
$\wh F(t^{-1})$, where $\wh F(t)=F^*(\bar t)$.
By applying (\ref{Ga}) and the above Lemma,
we deduce the following statement.

\begin{thm} 
Suppose $F$ be a matrix function in $\NDARNm$, where $m\ge1$. Let
$$
F^*(t)=h(t)\al^{-1}(t), \qquad t\in \BT
$$
be the right coprime factorization of $F^*$ on $\BT$,
where $\al$ is a Blaschke -- Potapov product
in $\BD$. Then
the space $M$ and matrix parameters $C$ and $\La$ of
the subnormal operator
$S=T_F$ (see (\ref{CLa})) can be calculated by the
formulas
\beqn
\label{MCLa}
M=\Htwom\ominus\al\Htwom, \quad
\La^*=T_{F^*}|M, \quad
C=\Ga_{F^*}^*\Ga_{F^*}|M.
\neqn
In particular, a pair $(C, \La)$ of operators on
a finite dimensional space $M$ gives rise to a
subnormal operator of finite type without point masses iff this
pair is unitarily equivalent to a pair $(C,\La)$, given by the above formulas.
\qed
\end{thm}

It is easy to write down an explicit orthonormal basis
of the space $\Htwom\ominus\al\Htwom$ (see
the Malmquist-Walsh lemma in \cite{Nikbook} and also
Example 1 below).
This permits one to calculate matrices $C$ and $\La$ explicitly.

Notice that if $F$
is a $m\times m$ rational matrix function in
$\Hinfmm$ with $\det F\not\equiv0$
and $F=\al h^*$ on the unit circle, where
$\al$, $h$ are right coprime and
$\al$ is a Blaschke--Potapov product,
then the symbol $F$ is normal (that is, $F^*F=FF^*$ on the
unit circle) if and only if $F=h\theta$, where $\theta$
is a $m\times m$ rational function, which is unitary on $\BT$,
but not necessarily analytic in $\BD$.

\section{A method of constructing rational matrix functions of classes
$\NDRNmm$ and $\NDARNmm$}

Let $B$ be a Blaschke--Potapov product, and let $\psi(t,\eta)$
be a scalar rational function of two variables. Put
\beqn
F(t)=\psi\big(t, B(t)\big),
\label{BlaschPot}
\neqn
and suppose that $F$ is well-defined
as a meromorphic function on the complex plane.
Since $B$ is unitary on the unit circle, $F$ is a
rational normal matrix function. For a fixed $B$ and ``most''
functions $\psi$, $F$ is non-degenerate, hence
a function of class
$\NDRNmm$. If, moreover, $F$ is analytic on the closed unit
disc, then $F$ belongs to $\NDARNmm$.
So (\ref{BlaschPot}) can be useful
in the construction of separated
real-type algebraic curves and quadrature domains.

In fact, it is more than an example: it can be proved that,
basically, any function in $\NDRNmm$ can be
obtained in the above way.
This topic will be treated in more detail elsewhere.

For functions $F$, obtained by this rule, the
Ahlfors type curve $\De_3(F)$ is closely related to the study of
the separated algebraic curve
\beqn
\label{rhoB}
\rho(B)=\{(t,\eta)\in\BC^2: \det\big(B(t)-\eta I\big)=0\}
\neqn
with the anti-analytic involution
\beqn
\label{invo2}
\de=(t,\eta)\mapsto \de^*=(\bar t^{-1},\bar\eta^{-1})
\neqn
(in general, it can be reducible).
The real part $\rho_{\BR}(B)$ of this curve is defined by the equation
$|t|=1$; it has a real dimension one.
The equality $|\eta|=1$ also holds true on $\rho_{\BR}(B)$.

Define the meromorphic function
$z=\psi(t,\eta)$ on $\rho(B)$. Then
the ``image'' $\ga(F)$ of the matrix function $F(t)$
on the unit circle $\BT$, which was defined by  (\ref{gammaF}),
equals to the $z$-image of the
real part $\rho_{\BR}(B)$ of the curve $\rho(B)$.
These notions will be exploited in
Example 2 of the next section.

\section{Concrete examples}

We put $b_\la(t)=\frac{t-\la}{\vph 1-\bar \la t}$.

\begin{exnum}
[matrix parameters of a simply connected quadrature domain]  
Take
$a,\be\in\BC$ with $|a|>1$ and $\be\ne0$, and put
$$
F(t)=t+\frac\be{t-a}, \qquad \Om=F(\BD).
$$
Then $F\in\NDARN_1$. Assume that $F$ is univalent on $\BD$ (for a
fixed $a$, it always can be achieved by taking a small $\be$). Then
$\Om$ is a quadrature domain by the Aharonov--Shapiro Theorem B. It corresponds to
the analytic Toeplitz operator $S=T_F$ on the scalar $H^2$, which, as we know,
is a subnormal operator of finite type without point masses. We are going to
calculate the matrix parameters $(C,\La)$ of this subnormal operator.

By Lemma \ref{LEM9} and Theorem 5, the space $M=\Range(S^*S-SS^*)$ can
be calculated as
$$
M=\big(\Ker \Ga_{\bar F}\big)^\perp.
$$
The function $\bar F$ coincides on $\BT$ with the rational function
$$
F_*(t)=\overline{F(\vphantom{\hat L}\bar t^{-1})}=t^{-1}+\frac {\bar\be t}{1-\bar a t}.
$$
This function has poles $0$ and $\bar a^{-1}$ in $\BD$ of order one.
It is easy to see that $\Ker \Ga_{\bar F}=\{x\in H^2: x(0)=x(\bar a^{-1})=0\}$.
Hence
$M=H^2\ominus \al H^2$, where $\al(t)=b_0(t)b_{\bar a^{-1}}(t)$, and
$\dim M=2$.

Next, put $M_-=\Range \Ga_{\bar F}\subset H^2_-$. Then $\dim M_-=2$. Choose
some orthonormal bases $\{e_1,e_2\}$ in $M$ and
$\{h_1,h_2\}$ in $M_-$. Then $P_-F_*e_j\in M_-$, $j=1,2$. Since
$S^*M\subset M$, $P_+F_*e_j\in M$ for $j=1,2$.
Hence there are expansions
$$
\label{expns}
\begin{aligned}
F_*e_1& = r_{11}h_1+ r_{21}h_2+\nu_{11}e_1+\nu_{21}e_2, \\
F_*e_2& = r_{12}h_1+ r_{22}h_2+\nu_{12}e_1+\nu_{22}e_2.
\end{aligned}
$$
Introduce an operator $R=\Ga_{F^*}\big|M: M\to M_-$, then
by Theorem 5, $C=R^*R$. By (\ref{MCLa}),
$R \sim
\begin{pmatrix}
r_{11} & r_{12} \\
r_{21} & r_{22}
\end{pmatrix}
$
and
$\La^* \sim
\begin{pmatrix}
\nu_{11} & \nu_{12} \\
\nu_{21} & \nu_{22}
\end{pmatrix}
$
in the bases $\{e_1,e_2\}$ and $\{h_1,h_2\}$.
In particular, one can take
$e_1=1$, $e_2={kt}(1-a^{-1}t)^{-1}$,
$h_1=t^{-1}$, $h_2=kt^{-1}(t-\bar a^{-1})^{-1}$, where
$k=\sqrt{1-|a|^{-2}}$. 
After calculating coefficients $r_{js}$,
$\nu_{js}$ (using residues), one gets
\beqn
\La^*=
\begin{pmatrix}
-\frac{\bar\be}{\bar a} \; & \;k-\frac{\bar\be}{\vvph\bar a^2 k} \\

\vph 0                  \;     & \;\frac 1 a-\frac{\bar \be}{\bar a k^{2} }
\end{pmatrix},
\quad
R=
\begin{pmatrix}
1-\frac {\bar\be}{\vvph\bar a^{2} }\; & \;  -\frac { \bar\be } {\vvph k\bar a^{3} } \\
\vph -\frac{ \bar\be } {\vvph k\bar a^{3} }\; & \;-\frac { \bar\be } {\vvph k^{2}\bar a^{4} }
\end{pmatrix}.
\neqn
Eigenvalues of $\La$ coincide with the nodes of the quadrature domain
$\Om$
(which are the points $F(0)$ and $F(\bar a^{-1})$).
It is also known from
the theorem by Helton--Howe and Carey--Pincus (see \cite{MP}) that
$$
\Area(\Om)=\pi\,\trace(C).
$$
In \cite{Ysubone}, we called $\La$ and $C^{1/2}$ the matrix
center and the matrix radius of $S$.

By Lemma 2, the operator of multiplication by
$F$ on $L^2(\BT)$ is the minimal normal extension
of the analytic Toeplitz operator $T_F$.
Since $\si(T_F)=F(\BT)$,  the curve
$\ga(F)=F(\BT)=\partial\Om$ can be described
alternatively by the equation (\ref{siN}).

Representations of the boundary of a quadrature domain by an equation
like (\ref{siN}) have been also considered in the series
of papers by Putinar and Gustafsson, see \cite{PUT}, \cite{GP} and
earlier papers. See also
Xia \cite{Xfour}, \cite{Xrev} and others.
These papers deal with hyponormal operators, instead
of subnormal ones. The domain need not be simply connected. The difference with
the subnormal case is that one can always find
a pair $(C,\La)$ that give rise to a quadrature domain so that
$\rank C=1$. In our setting, $C$
has always a full rank. The advantage of matrix parameters of
a quadrature domain in the sense of Putinar and Gustafsson
is that they are always determined uniquely. If a quadrature
domain is not simply connected, there are many subnormal operators
that correspond to it, and they give rise to different pairs
$(C,\La)$. The ambiguity is codified by the so-called characters
\cite{Fish}, \cite{Fed97}, \cite{Ysubtwo}.
\end{exnum}

\begin{exnum}[A one-connected quadrature domain] 
Consider a finite Blaschke product in $\BC^2$ of degree two:
$$
B(t)=\big( I-Q_1+b_\la(t)Q_1 \big)\big( I-Q_2+b_{-\la}(t)Q_2 \big),
$$
where $Q_1$, $Q_2$ are two different rank one projections in
$\BC^2$ and
$\la$ is a fixed point in $\BD$ with $\Ree \la\ne0$.
Assume (without loss of generality)
that $Q_j=\ell_j \otimes \ell_j$, where
$\ell_1=(1,0)$, $\ell_2=(c,a)$, with
$a> 0$, $c\ge0$, $a^2+c^2=1$.
Put
\begin{gather}
p(t)=1-\bar\la^2t^2, \qquad r(t)=t^2-\la^2, \\
q(t)= c^2(1-\bar\la^2)t^2+2a^2(1-|\la|^2)t+c^2(1-\la^2),
\end{gather}
and
\beqn
\label{psi}
\psi(t,\eta)=
\frac
{2p(t)\eta-q(t)+L(t-\ga_1)}
{2p(t)\eta-q(t)-L(t-\ga_1)}\, ,
\neqn
where $L\ne0$ is a complex constant, and
$\ga_1$ is a root of the polynomial $D(t)=q^2(t)-4p(t)r(t)$ with
$|\ga_1|<1$. We assert that there is a continuum
of
parameters $a,c,\la,\ga_1,L$ such that the matrix function
\beqn
\label{psiB}
F(t)=\psi(t,B(t))
\neqn
\begin{figure}
\includegraphics[width=140pt,height=200pt]{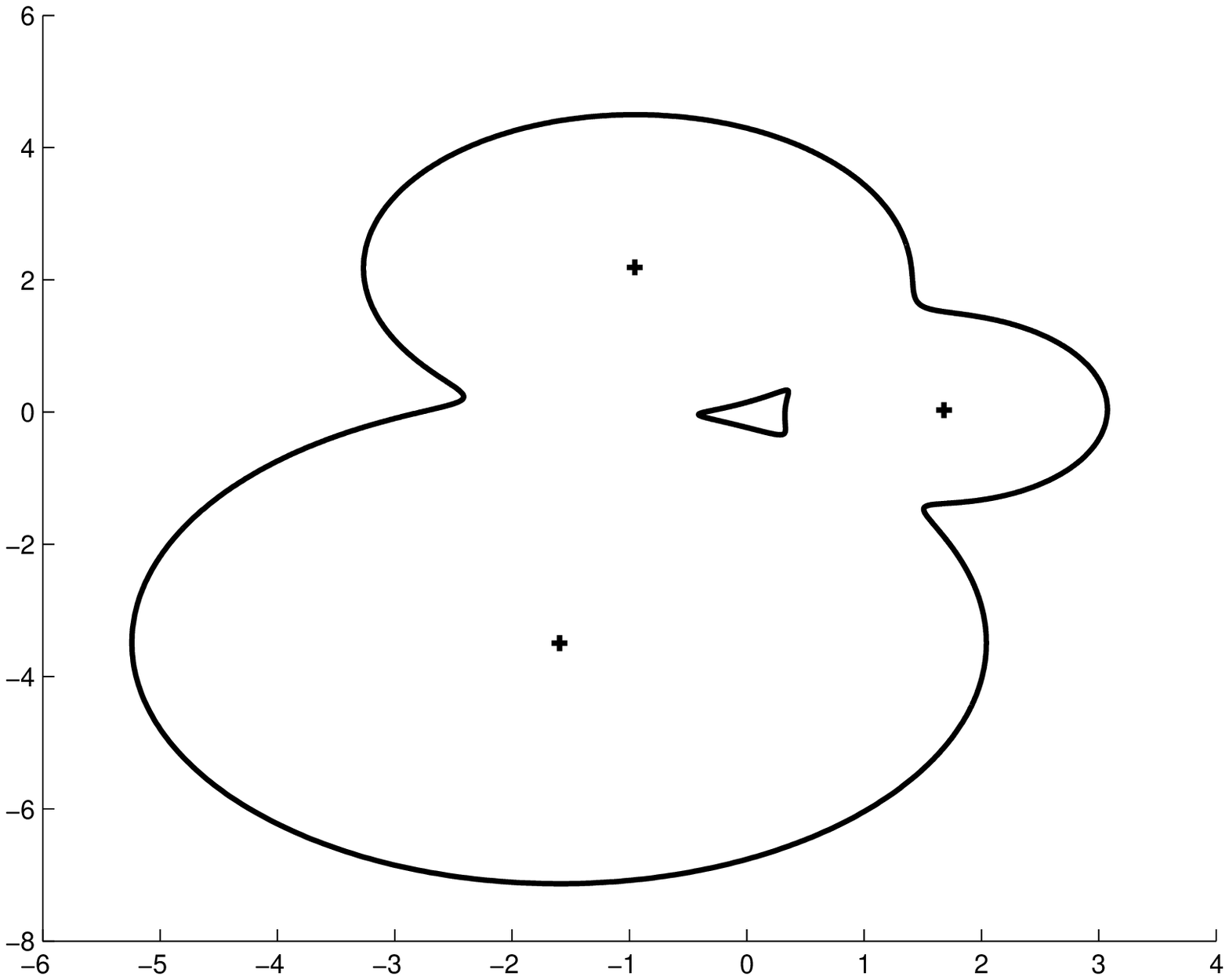}
\hskip2cm
\includegraphics[width=180pt,height=175pt]{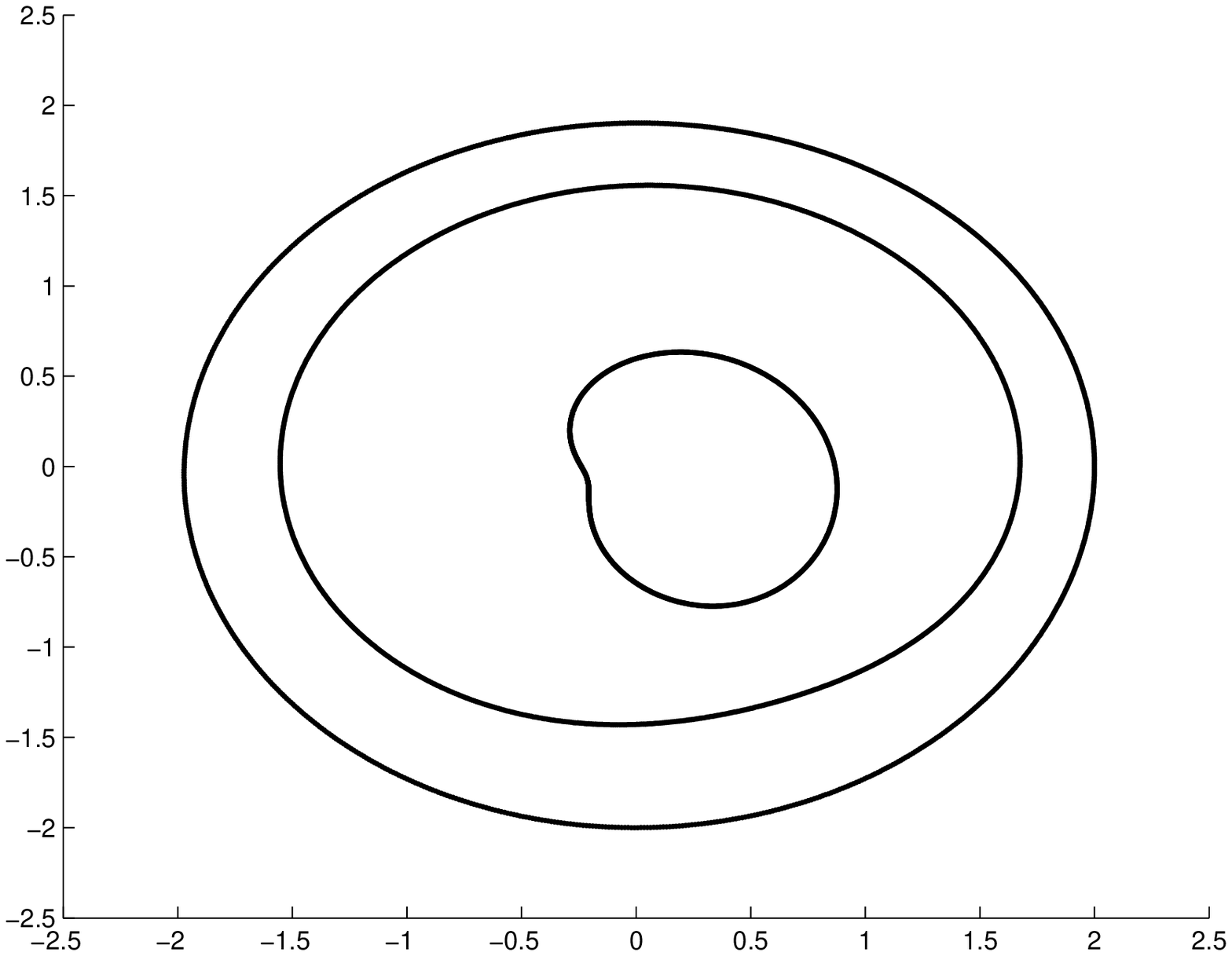}
\begin{center}
\hskip-2cm Fig. 1 \hskip6cm Fig. 2
\end{center}
\end{figure}
\noindent belongs to $\NDARN_2$ and gives rise to a
one-connected quadrature domain
according to the rule of Theorem 1.
It is so, in particular, if one chooses
$\la=0.8i$,
$a=\frac 5{13}$, $c=\frac {12}{13}$, $L=i$, and
the corresponding root $\ga_1 \thickapprox 0.0729 - 0.6467i$
of $D$.
Figure 1 shows the shape of the curve
$\ga(F)$ for these concrete parameters.
The explicit parametrization of the two parts of this curve is
$z=z_\pm(t)$, $t\in\BT$;
the functions
$z_\pm(t)$ will be defined in (\ref{B}),
(\ref{31}).
\end{exnum}

In what follows, we will motivate  this example and give more
comments and details about this quadrature domain.

A direct calculation shows that the algebraic curve (\ref{rhoB})
for our choice of $B$ has the form
\beqn
\label{curve}
\rho(B):\qquad p(t)\eta^2-q(t)\eta+r(t)=0.
\neqn
Obviously, $D(t)$ is the discriminant of this quadratic equation
in $\eta$. One observes that
$\overline{r(\vvph\bar t^{-1})}=t^{-2}p(t)$,
$\overline{q(\vvph\bar t^{-1})}=t^{-2}q(t)$, which
implies that $\overline{D(\vvph\bar t^{-1})}=t^{-4}D(t)$. Hence the
roots of $D$ are symmetric with respect to the unit circle. Denote
them as $\ga_1$, $\ga_2$, $\bar\ga_1^{-1}$, $\bar\ga_2^{-1}$. We
assume that all these roots are distinct and that $\ga_1, \ga_2\in
\BD$. The algebraic curve (\ref{curve}) is irreducible, and by
taking its normalization we can regard it as a compact Riemann
surface. Since the formula of the solution of (\ref{curve}) is
$\eta_{\pm}=\frac{q(t)\pm \sqrt{D(t)}}{2p(t)}$ and $p,q$ are
single-valued functions, the surface $\rho(B)$ coincides with the
Riemann surface of the multi-valued function
$$
t\mapsto \sqrt{D(t)}=K\cdot\sqrt{
(t-\ga_1)(t-\ga_2)(t-\bar\ga_1^{-1})(t-\bar\ga_2^{-1})}.
$$
Hence $\rho(B)$ is an elliptic curve and is homeomorphic to a
torus.

Inequalities $|t|<1$ and $|t|>1$ define the two
``halves'' of the curve $\rho(B)$, which we denote as $\rho_+(B)$ and
$\rho_-(B)$, respectively. Then $\rho_+(B)$ is homeomorphic to the
Riemann surface of the function $\sqrt{(t-\ga_1)(t-\ga_2)}$,
defined on the disc $|t|<1$. It is a two-sheeted branched covering
of the unit disc, and thus is homeomorphic to a sphere with
two holes (or to a one-connected domain).

Following B. Gustaffson \cite{Shdou}, pages 224--225,
we can search a meromorphic function $z(\cdot)$ on
$\rho(B)$, which has no poles on $\clos \rho_+(B)$ and is
univalent on $\rho_+(B)$. The image of  $\rho_+(B)$ under any such
function $z$ will be a one-connected quadrature domain.

Consider the meromorphic function
$z(\cdot)=\phi\big(t(\cdot),\eta(\cdot)\big)$ on $\rho(B)$. Note that
\beqn
\label{B}
z=\psi(t,\eta)=\frac{\si+L(t-\ga_1)}{\si-L(t-\ga_1)},
\neqn
where
\beqn
\label{31}
\si\defin 2p(t)\eta-q(t)=\pm\sqrt{D(t)}.
\neqn
By considering the local parameter  $\si$ on the curve
$\rho_+(B)$ in a neighborhood of the branching point
$t=\ga_1$, one gets that $z(\de)$ has no pole at this point.

There are two global continuous branches of
$\sqrt{D(t)}$ on the unit circle.
Define the functions $z_+(t)$, $z_-(t)$
on $\BT$ by putting in (\ref{B}) $\si=\pm\sqrt{D(t)}$, respectively.

The parameters $a,c,\la,\ga_1,L$
lead to a quadrature domain $\Om\defin z\big(\rho_+(B)\big)$
if and only if the following two conditions hold.
The first condition is that the
function $z(\de)$ should have no
poles on $\rho_+(B)$. The second one is that
$z(\cdot)$ should be univalent on $\rho_+(B)$
(or, equivalently, that
functions $\zeta_1(\tht)=z_-(e^{i\tht})$,
$\zeta_2(\tht)=z_+(e^{i\tht})$ should satisfy the
topological condition of Theorem 1). If these two conditions are valid,
then $\Om$ is a quadrature domain. In this case, we can agree that
when $t$ runs over the unit circle,
$z_+(t)$ traverses the outer boundary curve of $\Om$ and
$z_-(t)$ traverses the inner one. One gets from (\ref{B})
that $z_-(t)\cdot z_+(t)=1$  for $|t|=1$.

To verify the second condition for
a concrete set of values of parameters, it suffices to check,
for instance, that $\arg z_+(e^{i\tht})$ strictly increases for
$\tht \in[0,2\pi]$. The author has checked both
conditions numerically for the parameters indicated above.
It follows that close values of parameters also give a
quadrature domain.

In fact, the author does not know whether the
first necessary condition implies the second one.

Formulas (\ref{B}), (\ref{31})  were found
by the analogy with the inverse Zhukovsky function.

Since the meromorphic function $z=\psi(t,\eta)$ has no poles on
$\clos \rho_+(B)$, it follows that $F$ is analytic
on $|t|<1$. Hence $F\in\NDARN_2$.
Notice that for $t\in\BT$, $z_\pm(t)=\psi(t,\eta_\pm)$, where
$\eta_\pm$ are the two roots of the quadratic equation
(\ref{curve}). Hence $F(t)$ has eigenvalues
$z_+(t)$ and $z_-(t)$ for $|t|=1$.
It follows that, whenever our choice of parameters
produces a quadrature domain $\Om$, function $F$ generates the
same quadrature domain.
The spectrum of $T_F$ coincides with the closure of $\Om$.

\textbf{The Schwartz function and the nodes. The defining equation}
Suppose that our parameters
$\la, a, c, L, \ga_1$ are admissible, that is, they
give rise to a quadrature domain $\Om$. Then
the Schwartz function is given by
$$
w(z)=\overline {\psi(\de(z)^*)},\qquad z\in \clos\Om
$$
where $\de(z)$ is the function inverse to the function $z|\rho_+(B)$.
Function $w(z)$ has three poles, which
are the nodes of the quadrature domain (the points $z_j$ in the
formula (\ref{tAA})). The positions of these three nodes are indicated
on Fig. 1.

It is possible (in principle)
to write down explicitly the polynomial defining equation
of this quadrature domain. Namely, one can derive a polynomial relation
$
X_z(t)\equiv \sum_{j=1}^K X_j(z)t^j=0
$
between the functions $t$, $z$ on the curve $\rho(B)$
(here $X_j$ are polynomials of one variable).
Next, one can write down a similar polynomial relation
$
Y_w(t)\equiv \sum_{j=1}^K Y_j(w)t^j=0
$
between the functions $t$, $w$ on $\rho(B)$.
Then one has an explicit equation
$$
\Res(X_z,Y_w)=0,
$$
which satisfy the meromorphic ``coordinates''  $z$, $w$
on $\rho(B)$ and which is polynomial in $z$, $w$. Here
$\Res(M,N)$ is the resultant of polynomials $M(t)$, $N(t)$, see
\cite{Waerd}. We recall that
$\Res(M,N)$ vanishes iff $M$ and $N$ have a common root. It is not
completely clear, however, whether this equation is
a minimal one (it might have extra factors of the form $z-z_0$ or $w-w_0$).

Numerical experiments show that for different values of
admissible parameters, the quadrature domain obtained has a form of
three merged circular drops, and the nodes of the domain are situated
approximately in the centers of these drops.
This justifies the following
\begin{conj}
If the above method leads to a quadrature domain, then this domain
can be obtained alternatively as a final domain
at a time $T_0>0$ from a Hele--Shaw flow
with three sources, situated in the three nodes of the domain
(with no liquid at the starting time $T=0$).
\end{conj}
We refer to the book \cite{Sak} for a discussion of the Hele--Shaw
flows and to
\cite{VarchEtin} for a good elementary introduction to the subject.
In \cite{Cro3}--\cite{Cro}, \cite{Richsn},
one can find more
recent results.

A general characterization of one-connected
quadrature domains was given by M. Putinar in
\cite{PUTJFA}, Theorem 1.3. However, he
did not give concrete examples.

\begin{exnum}[A matrix function of class $\NDARN_3$]
One can try to use the same ideas in order to obtain
more complicated examples. For any $m$, it is easy to
get many examples of functions in $\NDARN_m$. For instance,
put $Q_j=\ell_j\otimes \ell_j$, where
$\ell_1=(\frac  {12}{13}, \frac {-5}{13}, 0)$,
$\ell_2=(0, \frac {12}{13}, \frac {-5}{13} )$,
$\ell_3=(\frac {-5}{13}, 0, \frac {12}{13})$. Put
$\eps_1=i$ and $\eps_2=\exp(\frac 23\pi i)$.
Let $\la=\frac 1{10}\in\BD$.
Take the matrix Blaschke product
$$
B(t)=\big( I-Q_1+\eps_1b_\la(t)Q_1 \big)
\big( I-Q_2+\eps_2b_{-\la}(t)Q_2 \big)
\big( I-Q_3+tQ_3 \big)
$$
and the function $\psi(t,\eta)=t+\eta$. It is easy to see that
the matrix function $F(t)=\psi(t,B(t))$ is in $\NDARN_3$.
The corresponding curve $\ga(F)$ is shown on Fig. 2.
It can be deduces from this picture that the real part of
the algebraic curve $\rho(B)$ has three components.
One can number the three eigenvalues of $F(e^{i\tht})$
so that when $\tht$ runs over $[0,2\pi]$, each of these eigenvalues
traverses its own component of the curve $\ga(F)$ in the
positive direction.
This function $B$ can be thought of as a perturbation
of the case when $\{\ell_j\}$ are an orthonormal basis. In the latter case,
$\rho(B)$ has three irreducible pieces, and $\ga(F)$ consists of
three concentric circles.

If we had taken
$\eps_2=\frac {-1+i}{\sqrt{2}}$,
without changing other data, then
$\rho_{\BR}(F)$ would have only two components.
\end{exnum}
In general, it is not so easy to
find out the topological types of curves
$\rho(B)$. If $\rho_+(B)$ is homeomorphic to
a multiply connected domain, one could try to
construct meromorphic functions $z=\psi(t,\eta)$ on $\rho(B)$, which
give rise to quadrature domains.
It is unclear by now how to do it explicitly.
It would be desirable to have some general results
about possible topological types of $\rho(B)$, depending on
the size and the degree of a Blaschke--Potapov product
$B$.

A general method of constructing multiply connected quadrature
domains has been suggested recently by Crowdy \cite{Cro2}, \cite{Cro3},
by Crowdy and Marshall in
\cite{Cro} and by Richardson in \cite{Richsn}.
In particular, the work by Crowdy and Marshall
contains many examples of calculation of
quadrature domains of different connectivity.
The method by Crowdy involves Schottky--Klein prime functions,
defined as an infinite product.
Richardson's method uses Poincar\'e series. Both
methods require to solve certain systems of
nonlinear equations. It would be interesting to find relations
between the methods by Crowdy, Marshall and Richardson and the results of
the present article.

\section{Further perspectives}
One can try to apply our results to several neighboring fields.
There are many unanswered questions about quadrature domains, see
the collection of problems in the recent book \cite{EbenGuKhPutinar}.
Our Theorem 1, combined probably with some new ideas, could
occasion some progress. Let us indicate two concrete problems in
this connection.

1)
As it is proven in \cite{Shdou}, Thm. 12, for any $p\ge 1$,
there is a family of $p$-connected domains that satisfy the same quadrature identity
(\ref{tAA}) and depend on at least $p$ real parameters.
It would be
interesting to find a more or less explicit parametrization of all
quadrature domains of a given connectivity
satisfying the same quadrature identity.

It is not known whether there is uniqueness
when one considers only simply connected domains; see
\cite{Sakai3} for a partial result.

2) A point $z_0$ of a quadrature domain $\Om$
is called \textit{special} if
$w(z_0)=\bar z_0$, where $w(z)$ is the Schwartz function.
If $\Om$ is generated by a matrix function $F$ of class $\NDARN$, then
$z_0$ has to be an eigenvalue of $F(t_0)$ and
$\bar z_0$ an eigenvalue of $F^*(\bar t_0^{-1})$ for a point
$t_0$ in $\BD$. It is interesting to estimate the number of special points
of a quadrature domain; see
\cite{Gust-sing} and \cite{Sakai4} for some results in this direction.
We remark that special points play an important role in the
connection between subnormal and hyponormal operators attached to
$\Om$, see \cite{Yqua}, Thms. 1, 2 and the proof of Lemma 3.

The connection between hyponormal operators and quadrature domains
has been exploited in works by Gustafsson,
Putinar and Xia, see
\cite{PUTJFA}, \cite{PUT}, \cite{GP}, \cite{Xrev},
\cite{Yqua} and references therein.

Certainly, a better understanding of ways to construct
multiply connected quadrature domains would be
very desirable. A concrete construction of Ahlfors functions
on a multiply connected domain has been given by S. Fedorov
in \cite{Fed90}.

It is also interesting to look for numerical applications of
quadrature domains.
They depend on a finite number of parameters and
can approximate an arbitrary bounded domain,
see \cite{GustHeMiPutinar} for a discussion.
Our results might help in constructing conceptually
simple algorithms of dealing with quadrature domain.

There are many topics that are related to the subject of this
work, which we did not touch.
Hyponormal Toeplitz operators with finite rank self-commutator
were studied in \cite{Cow}, \cite{HwKLee} and other papers;
a relationship between Toeplitz operators with rational symbol and Riemann surfaces
was exploited in \cite{YRSm}. Multiplication operators by the independent
variable and by an analytic matrix function were studied in
\cite{Voich}, \cite{Hasu}, \cite{AbDo}, \cite{VS1}, \cite{VS2},
\cite{YRSm} and others.

Vector bundles over real algebraic curves were used here very little.
The topic of Section 7, in fact, is related to the
so-called determinantal representations of vector bundles
of real algebraic curves and characters. See \cite{AbDo}, \cite{Vinn1},
\cite{Vinn2} and the references therein.

Real algebraic curves and vector bundles over them
appear naturally in the theory of commuting nonselfadjoint operators
and discrete and continuous linear systems with
multidimensional time, which is being developed by
Liv\v sic, Alpay, Ball, Vinnikov and others.
We refer to \cite{AlpV}, \cite{BallVin} (and the references therein) and to the book
\cite{LKMV}. As it was shown in \cite{AlpV},
there are advantages in defining spaces $H^2$ as spaces
of differentials of order $1/2$ instead of spaces of functions.
Some algebraic and
computational aspects of this theory
were developed in \cite{AShap}.

Curves as in (\ref{rhoB}), whose anti-analytic involution
has the form
$(t,\mu)\mapsto (\bar t^{-1},\bar\mu^{-1})$
appear in the theory of commuting contractions, see
\cite{BallVin2}.

In a series of works,  Pavlov and Fedorov developed
the harmonic analysis on multiply connected domains
(see \cite{{Fed90}}, \cite{{Fed97}} and others). Such
topics as analogues of Muckenhoupt condition and of Carleson
condition, Ahlfors type functions
that generate the uniform analytic algebra
in the domain,
coinvariant subspaces and corresponding semigroups
were studied by these authors. One of the aims of this program is
to develop a kind of the Lax-Phillips approach to the investigation
of resonances for a selfadjoint operator with band spectrum, see
\cite{{Pav}}.

Algebraic curves also appear systematically in the study of
integrable dynamical system, which is a very vast
area; we only mention the review \cite{DubKrN}.
As the work \cite{LivAv} suggests, this topic also has
strong connections with the theory of
commuting nonselfadjoint linear operators.



\end{document}